\documentclass[12pt]{article}
\usepackage[left=1in,right=1in,top=1in,bottom=1in]{geometry}
\usepackage{cite}
\usepackage{graphicx,bbm,pstricks,soul}
\usepackage{pifont}
\usepackage{bbm,algorithmic,mdframed,placeins,multirow,booktabs,subfigure}
\usepackage[ruled,vlined,linesnumbered]{algorithm2e}

\usepackage{amsmath,amsfonts,amsbsy,amssymb}

\newcommand{\LRARR}[4]{{\mbox{ \raise 0.4 mm \hbox{$#1$}}} \;
  \mathop{\stackrel{\displaystyle\longrightarrow}\longleftarrow}^{#3}_{#4}
  \; {\mbox{\raise 0.4 mm\hbox{$#2$}}}}

\newcommand{\G}{\mathcal{G}}

\newcommand{\X}{{\mathbf X}}
\newcommand{\Y}{{\mathbf Y}}
\newcommand{\W}{{\mathbf W}}
\newcommand{\data}{D}
\newcommand{\neff}{n_{\text{eff}}}
\newcommand{\E}{{\mathbb E}}
\renewcommand{\b}[1]{{\bf #1}}
\DeclareMathOperator*{\argmin}{arg\,min}
\DeclareMathOperator*{\argmax}{arg\,max}

\newcommand{\picturesAB}[6]{
\centerline{
\hskip #4
\raise #3 \hbox{\raise 0.9mm \hbox{(a)}}
\hskip #5
\epsfig{file=#1,height=#3}
\hskip #6
\raise #3 \hbox{\raise 0.9mm \hbox{(b)}}
\hskip #5
\epsfig{file=#2,height=#3}
}}
\newcommand{\picturesCD}[6]{
\centerline{
\hskip #4
\raise #3 \hbox{\raise 0.9mm \hbox{(c)}}
\hskip #5
\epsfig{file=#1,height=#3}
\hskip #6
\raise #3 \hbox{\raise 0.9mm \hbox{(d)}}
\hskip #5
\epsfig{file=#2,height=#3}
}}

\makeatletter  
\newcommand{\xleftrightarrows}[2][]{\mathrel{%  
 \raise.40ex\hbox{$  
       \ext@arrow 3095\leftarrowfill@{\phantom{#1}}{#2}$}%  
 \setbox0=\hbox{$\ext@arrow 0359\rightarrowfill@{#1}{\phantom{#2}}$}%  
 \kern-\wd0 \lower.4ex\box0}}  
 
\newcommand{\xrightleftarrows}[2][]{\mathrel{%  
 \raise.40ex\hbox{$\ext@arrow 3095\rightarrowfill@{\phantom{#1}}{#2}$}%  
 \setbox0=\hbox{$\ext@arrow 0359\leftarrowfill@{#1}{\phantom{#2}}$}%  
 \kern-\wd0 \lower.4ex\box0}}  
 
\def\leftrightarrowfill@{%
 \arrowfill@\leftarrow\relbar\rightarrow%
 }
\newcommand*{\centerfloat}{%
  \parindent \z@
  \leftskip \z@ \@plus 1fil \@minus \textwidth
  \rightskip\leftskip
  \parfillskip \z@skip}
\makeatother
\makeatother 

\author{Colin Cotter\thanks{Department of Mathematics, Imperial
    College, London, UK} \and Simon Cotter\thanks{School of
    Mathematics, University of Manchester, Manchester, UK. e:
    simon.cotter@manchester.ac.uk. SLC is grateful for EPSRC First
    grant award EP/L023393/1} \and Paul Russell\thanks{School of
    Mathematics, University of Manchester, Manchester, UK}}
\title{Ensemble Transport Adaptive Importance Sampling}
\begin{document}
\maketitle
\begin{abstract}
Markov chain Monte Carlo methods are a powerful and commonly used
  family of numerical methods for sampling from complex probability
  distributions. As applications of these methods increase in size and
  complexity, the need for efficient methods increases. In this
paper, we present a particle ensemble algorithm. At each iteration, an
  importance sampling proposal distribution is formed using an
  ensemble of particles. A stratified sample is taken from this
  distribution and weighted under the posterior, a state-of-the-art
  ensemble transport resampling method is then used to create an evenly weighted sample
  ready for the next iteration. We demonstrate that this ensemble transport
  adaptive importance sampling (ETAIS) method outperforms MCMC methods
  with equivalent proposal distributions for low dimensional problems,
  and in fact shows
  better than linear improvements in convergence rates with respect to
  the number of ensemble members. We also introduce a new resampling
  strategy, multinomial transformation (MT), which while not
  as accurate as the ensemble transport resampler, is substantially less costly for large
  ensemble sizes, and can then be used in conjunction with ETAIS for
  complex problems. We also focus on how algorithmic parameters
  regarding the mixture proposal can be quickly tuned to optimise
  performance. In particular, we demonstrate this methodology's
  superior sampling for multimodal problems, such as those arising
  from inference for mixture models, and for problems with expensive
  likelihoods requiring the solution of a differential equation, 
    for
  which speed-ups of orders of magnitude are demonstrated. Likelihood
  evaluations of the ensemble could be computed in a distributed
  manner, suggesting that this methodology is a good candidate
  for parallel Bayesian computations.
\end{abstract}
\section{Introduction}
Having first been developed in the early 1970s\cite{hastings1970monte}, Markov chain Monte Carlo (MCMC) methods have been of increasing
importance and interest in the last 20 years or so. They allow us to
sample from complex probability distributions which we would not be
able to sample from directly. In particular, these methods have
revolutionised the way in which inverse problems can be tackled,
allowing full posterior sampling when using a Bayesian framework. 

However, this often comes at a very high cost, with a very large
number of iterations required in order for the empirical approximation
of the posterior to be considered good enough. As the cost of
computing likelihoods can be extremely large, this means that many
problems of interest are simply computationally intractable.

This problem has been tackled in a variety of different ways. One
approach is to construct increasingly complex MCMC methods which are
able to use the structure of the posterior to make more intelligent
proposals, leading to more thorough exploration of the posterior with
fewer iterations. For example, the Hamiltonian or Hybrid Monte Carlo
(HMC) algorithm uses gradient information and symplectic integrators
in order to make very large moves in state with relatively high
acceptance probabilities\cite{sexton1992hamiltonian}. Non-reversible
methods are also becoming quite popular as they can improve
mixing\cite{bierkens2015non}. Riemann
manifold Monte Carlo methods exploit the Riemann geometry of the
parameter space, and are able to take advantage of the local structure
of the target density to produce more efficient MCMC
proposals\cite{girolami2011riemann}. This methodology has been
successfully applied to MALA-type proposals and methods which exploit
even higher order gradient information\cite{bui2014solving}.

Since the clock speed of an individual processor is no longer
following Moore's law\cite{moore1998cramming}, improvements in
computational power are largely coming from the parallelisation of
multiple cores. As such, the area of parallel MCMC methods is becoming
increasingly of interest. One class of parallel MCMC method uses
multiple proposals, with only one of these proposals being
accepted. Examples of this approach include multiple try
MCMC\cite{liu2000multiple} and ensemble MCMC\cite{neal2011mcmc}. In
\cite{calderhead2014general}, a general construction for the
parallelisation of MCMC methods was presented, which demonstrated
speed ups of up to two orders of magnitude when compared with serial
methods on a single core. Another approach involves pre-fetching, where the possible future acceptances/rejections are calculated in advance\cite{angelino2014accelerating}.

A variety of other methods have been designed with particular scenarios in mind. For
instance, sampling from high/infinite-dimensional
posterior distributions is of interest in many applications. The majority of
Metropolis-Hastings algorithms suffer from the curse of
dimensionality, requiring more samples for a given degree of accuracy
as the parameter dimension is increased. However some
dimension-independent methods have been developed, based on
Crank-Nicolson discretisations of certain stochastic differential
equations\cite{cotter2013mcmc}. Other ideas such as chain
adaptation\cite{haario2005componentwise} and early rejection of
samples can also aid reduction of the computational
workload\cite{solonen2012efficient}.

However, high dimensionality is not the only challenge that we may
face. Complex structure in low
dimensions can cause significant issues. These issues may arise due to
large correlations between parameters in the posterior, leading to
long thin curved structures which many standard methods can struggle
with. These features are common, for example, in inverse problems
related to epidemiology and
other biological applications\cite{house2016bayesian}. 

Multimodality of the posterior can also lead to
incredibly slow convergence in many methods. Many methods allow for
good exploration of the current mode, but the waiting time to the next
switch of the chain to another mode may be large. Since many switches
are required in order for the correct weighting to be given to each
mode, and for all of the modes to be explored fully, this presents a
significant challenge. One application where this is an ever-present problem is that of
mixture models. Given a dataset, where we know that the data is from
two or more different distributions, we wish to be able to identify
the parameters, e.g. the mean and variance and relative weight, of
each part of the mixture\cite{marin2005bayesian}. The resulting posterior
distribution is invariably a multimodal distribution, since the likelihood is
invariant to permutations. Metropolis-Hastings algorithms, for
example, will often fail to converge in a reasonable time frame for
problems such as this. Since the posterior may be multimodal,
independent of this label switching, it is important to be able
to efficiently sample from the whole posterior.

Importance samplers are another class of methods which allow
us to sample from complex probability distributions. A related class of algorithms, adaptive importance sampling (AIS)
\cite{liu2008monte} reviewed in\cite{bugallo2015adaptive}, had
received less attention until their practical applicability was
demonstrated in the
mid-2000s\cite{celeux2006iterated,cappepopulation,isard1998condensation,bink2008bayesian}.
AIS methods produce a sequence of approximating distributions,
constructed from mixtures of standard distributions, from which samples can be
easily drawn. At each iteration the samples are weighted, often using
standard importance sampling methods. The weighted samples are used to
train the adapting sequence of distributions so that samples are drawn
more efficiently as the iterations progress. The weighted samples form
a sample from the posterior distribution under some mild
conditions~\cite{robert2013monte,martino2015adaptive}.

Ensemble importance sampling schemes also exist, e.g. population Monte Carlo (PMC)~\cite{cappe2012population}. PMC uses an ensemble to build a mixture or kernel density estimate (KDE) of the
posterior distribution. The efficiency of
this optimisation is restricted by the component kernel(s) chosen, and
the quality of the current sample from the posterior\cite{cappe2008adaptive,douc2007convergence,douc2007minimum}. Extensions, such
as layered adaptive importance sampling\cite{martino2017layered}, adaptive multiple importance sampling algorithm
(AMIS)~\cite{cornuet2012adaptive}, and adaptive population importance sampling
(APIS)~\cite{martino2015adaptive} have enabled these methods to be
applied to various applications, including population
genetics~\cite{siren2011reconstructing}.

In this paper, we present a framework for ensemble
importance sampling, which can be built around many of the current
Metropolis-based methodologies in order to create an efficient target
distribution from the current
ensemble. The method makes use of a resampler based on optimal
transport which has been used in the context of particle
filters\cite{reich2013nonparametric}. We also detail how algorithmic
parameters can be quickly tuned into optimal regimes, with respect to
the effective sample size statistic. In particular we demonstrate the
advantages of this method when attempting to sample from multimodal
posterior distributions, such as those arising from inference for
mixture models. The method that we present is specifically designed to
tackle the challenges of certain types of low-dimensional inverse
problem which have complex posterior structure, for example
multimodality or strong correlations, and expensive
likelihood evaluations. We will demonstrate that this methodology
allows us to sample from such a posterior distribution with greater
precision per likelihood evaluation than 
standard MCMC chains of the same type.

The method considered is similar to population Monte Carlo methods
which have been previously studied, for example in
\cite{elvira2017improving}. However our approach exploits the current
state-of-the-art resamplers, and adaptively optimises algorithmic parameters. Higher accuracy resamplers, that provide
a deterministic resample which optimally reflects the statistics of
the original sample, lead to better mixture approximations of the
posterior. These approximations can then be used to construct stable
and efficient importance sampling proposals. The cost of the
state-of-art resamplers can be prohibitive for large ensemble sizes,
which are necessary for sampling in moderately higher dimensions, and
as such we also present a greedy approximation of the optimal
transport resampler, which provides us with good results for a
fraction of the cost. We also detail how scaling parameters in the
MCMC proposals can be quickly tuned, leading us to an efficient and
fully-automated algorithm. We will demonstrate that despite additional
overheads, the ETAIS methodology can outperform its plain
Metropolis-Hastings cousin by orders of magnitude to reach a histogram
of a particular degree of convergence.

In Section \ref{Sec:Prelim} we introduce some mathematical
preliminaries upon which we will later rely. In Section \ref{Sec:ETAIS}
we present the general framework of the ETAIS algorithm. In Section
\ref{Sec:adapt} we consider adaptive versions of ETAIS which
automatically tune algorithmic parameters concerned with the proposal
distributions. In Section \ref{sec:MT} we introduce the multinomial transformation (MT) algorithm which is a less accurate but
faster alternative to resamplers which solve the optimal transport
problem exactly. In Section~\ref{sec:consistency} we consider
consistency of the ETAIS algorithm. In Section~\ref{sec:useful} we briefly look at one
particular advantageous property of this approach. In Section \ref{Sec:Num} we present
some numerical examples, before a brief conclusion and discussion in
Section \ref{Sec:Conc}.

%%%%%%

\section{Preliminaries}\label{Sec:Prelim}

In this Section we will introduce preliminary topics and algorithms
that will be referred to throughout the paper.

%%%%%%

\subsection{Bayesian inverse problems}
In this paper, we focus on the use of MCMC methods for characterising
posterior probability distributions arising from Bayesian inverse problems. We
wish to learn about a particular unknown quantity $x$, of which we are
able to make direct or indirect noisy observations. For now
we say that $x$ is a member of a Hilbert
space $X$. 

The quantity $x$ is mapped on to observable space by the observation
operator $\mathcal{G}:X \to\mathbb{R}^d$. We assume that the
observations, $\data$, are subject to Gaussian noise,
\begin{equation}\label{eqn:obs}
	\data = \mathcal{G}(x) + \varepsilon, \qquad \varepsilon \sim \mathcal{N}(0,\Sigma).
\end{equation}

These modelling assumptions allow us to construct the 
likelihood of observing the data $\data$ given the quantity $x =
x^*$. Rearranging \eqref{eqn:obs} and using the distribution of
$\varepsilon$, we get:
\begin{equation}\label{eqn:like}
	\mathbb{P}(\data|x=x^*) \propto \exp \left ( -\frac{1}{2} \|\mathcal{G}(x^*)
	  - \data\|_\Sigma^2 \right ) = \exp\left(-\Phi(x^*)\right),
\end{equation}
where $\| y_1 - y_2 \|^2_\Sigma = (y_1-y_2)^\top\Sigma^{-1}(y_1-y_2)$
for $y_1,y_2 \in \mathbb{R}^d$.

Then, according to Bayes' theorem, the posterior density that we
are interested in characterising is given, up to a constant of proportionality, by
\begin{equation*}
	\pi \propto \exp \left ( -\Phi(x^*) \right )\pi_0(x^*),
\end{equation*}
where $\pi_0$ is the prior density on the quantity of interest $x^*$.

% As discussed in \cite{stuart2010inverse,cotter2009bayesian},
% in order for this inverse problem to be well-posed in the Bayesian
% sense, we require the posterior distribution, $\mu_Y$, to be absolutely
% continuous with respect to the prior, $\mu_0$. A
% minimal regularity prior can be chosen informed by regularity results
% of the observational operator $\mathcal{G}$. Given such a prior, then
% the Radon-Nikodym derivative of the posterior measure, $\mu_Y$, with
% respect to the prior measure, $\mu_0$, is proportional to the
% likelihood:
% \begin{equation}\label{eqn:RND}
% 	\frac{d\mu_Y}{d\mu_0} \propto \exp \left ( -\Phi(x^*) \right ).
% \end{equation}

% %%%%%%

\subsection{Particle filters and resamplers}\label{sec:filters}
In this subsection we briefly review particle filters, since the
development of the resampler that we incorporate into the ETAIS is
motivated by this area. Particle filters are a class of Monte Carlo algorithm designed to
solve the filtering problem. That is, to find the best estimate of the
true state of a system when given only noisy observations of the
system. The solution of this problem has been of importance since the
middle of the 20th century in fields such as molecular biology,
computational physics and signal processing. In recent years the data
assimilation community has contributed several efficient particle
filters, including the ensemble Kalman filter
(EnKF)~\cite{evensen1994sequential} and the ensemble transform
particle filter (ETPF)~\cite{reich2013nonparametric}.

The ETPF defines a coupling $T$ between two random variables $Y$ and
$X$, allowing us to use the induced map as a resampler. An \emph{optimal
coupling} $T^*$ is one which maximises the correlation between $X$ and
$Y$ \cite{cotter2012ensemble}. This coupling is the solution to a
linear programming problem in $M^2$ variables with $2M-1$ constraints,
where $M$ is the size of the sample. Maximising the correlation
preserves the statistics of $X$ in the new sample.

In this work we use the resampler used within these particle filters. We
approximate the posterior density $\pi$, with an ensemble
of weighted samples, $\{(w_i,y_i)\}_{i=1}^M$. In filtering problems
the weights would be found by incorporating new observed data whereas
here we simply use importance weights. These weighted samples can then be resampled, using the ETPF or otherwise,
into a set of new equally-weighted samples $\{x_i\}_{i=1}^M$, such that
\begin{equation}\label{eqn:resampler}
	\sum\limits_{i=1}^M \! w_i\delta_{y_i}(\cdot)
	 \quad \xrightarrow{\text{ETPF}} \quad \sum\limits_{i=1}^M \!
	 \frac{1}{M} \delta_{x_i}(\cdot),
\end{equation}
where $\delta_{x}(\cdot)$ is the Dirac delta measure.

The resamplers used within particle filters such as the ETPF are well suited to this problem
since it is easy to introduce conditions in the resampling to ensure
you obtain the behaviour you require. One downside is that the
required ensemble size increases quickly with dimension, making it
difficult to use in high-dimensional problems.

%%%%%%

\subsection{Deficiencies of Metropolis-type MCMC schemes}
All MCMC methods are naively parallelisable. One can take a method
and simply implement it multiple times over an ensemble of processors. All
of the states of all of the ensemble member can be recorded, and in the time
that it takes one MCMC chain to draw $N$ samples, $M$ ensemble members
can draw $NM$ samples. 

However, we argue that this is not an optimal
scenario. First of all, unless we have a lot of information about the
posterior, we will initialise the algorithm's initial state in the
tails of the distribution. The samples that are initially made as the
algorithm finds its way to the mode(s) of the distribution cannot be
considered to be samples from the target distribution, and must be thrown
away. This process is known as the burn-in. In a naively parallelised
scenario, each ensemble member must perform this process independently,
and therefore mass parallelisation makes no inroads to cutting this cost.

Moreover, many MCMC algorithms suffer from poor mixing, especially in
multimodal systems. The number of samples that it takes for an MCMC
trajectory to switch between modes can be large, and given that a large
number of switches are required before we have a good idea of the relative
probability densities of these different regions, it can be prohibitively
expensive.

Another aspect of Metropolis-type samplers is that information
computed about a proposed state is simply lost if we choose to reject
that proposal in the Metropolis step. An advantage of importance
samplers is that no evaluations of $\mathcal{G}$ are ever wasted
since all samples are saved along with their relative weighting.

These deficiencies of MCMC methods
motivated the development of the Ensemble Transport Adaptive Importance Sampler
(ETAIS).
In the next section we will introduce the method in its most
general form.

%%%%%%

\section{The Ensemble Transform Adaptive Importance Sampler \allowbreak (ETAIS)}\label{Sec:ETAIS}

Importance sampling can be a very efficient method for sampling from a
probability distribution. A proposal density is chosen, from which we
can draw samples. Each sample is assigned a weight given by the
ratio of the target density and the proposal density at that
point. They are efficient when the proposal density is concentrated in
similar areas to the target density, and incredibly inefficient when
this is not the case. The aim of the ETAIS is to use an ensemble of states to construct a proposal
distribution which will be as close as possible to the target density. If this
ensemble is large enough, the distribution of states will be
approximately representative
of the target density.

The proposal distribution could be constructed in many different ways,
but we choose to use a mixture distribution, made up of a sum of MCMC
proposal kernels  e.g. Gaussians centred at the current state, as in RWMH. Once the proposal is constructed, we can
sample a new ensemble from the proposal distribution, and each is
assigned a weight given by the ratio of the target density and the
proposal mixture density. Assuming that our proposal distribution is a
good one, then the variance of the weights will be small, and we will
have many useful samples. Finally, we need to create a set of evenly
weighted samples which best represent this set of weighted samples.
This is achieved by implementing a resampling algorithm. These samples
are not stored in order to characterise the posterior density, since
the resampling process is not exact. They are
simply needed in order to inform the mixture proposal distribution for the next
iteration of the algorithm. 

Initially we
will use the ETPF resampler algorithm\cite{reich2013nonparametric}, although we
will suggest an alternative strategy in Section \ref{sec:MT}, for
examples where a
large ensemble is required, and for which the ETPF may become more expensive. The
resampling algorithm gives us a set of evenly weighted samples which
represents the target distribution well, which we can use to iterate
the process again. The algorithm is summarised in Algorithm
\ref{alg:ETAIS}.

The choice of resampling algorithm is important since its
  output is used to formulate the proposal distribution for the next
  iteration. A basic resampler may be cheap to implement, but the
  resulting sample, and in turn the next iteration's proposal distribution, may not be as representative of the target
  distribution. This leads to higher variances of the weights in the
  importance sampler, and therefore to slower convergence to the
  target density.

We wish to sample states $x \in X$ from a posterior
probability distribution with density $\pi$. Since we have $M$ ensemble members, we
represent the current state of all of the Markov chains as a vector
$\X = [x_1,x_2,\ldots,x_M]^\top$. We are also given a transition kernel
$\nu(\cdot,\cdot)$, which might come from an MCMC method, for example
the random walk Metropolis-Hastings proposal density $\nu(\cdot,x) \sim
\mathcal{N}(x,\beta^2)$, where $\beta^2\in \mathbb{R}$ defines the
variance of the proposal.

\begin{table}[!ht]
\centering
\begin{algorithm}[H]
\DontPrintSemicolon
\BlankLine
	Set $\X^{(0)} = \X_0 = [x_1^{(0)},x_2^{(0)},\ldots,x_M^{(0)}]^\top$.\;
	\For {$i = 1, \dots, N$}
	{
		Sample $\Y^{(i)} = [y_1^{(i)},y_2^{(i)},\ldots,y_M^{(i)}]^\top, \quad y_j^{(i)} \sim
\nu(\cdot;x_j^{(i-1)})$.\label{algline:ETAIS_propose}\;
		Calculate $\W^{(i)} = [w_1^{(i)},w_2^{(i)},\ldots,w_M^{(i)}]^\top,$ \quad $w^{(i)}_j =
\frac{\pi(y_j^{(i)})}{\chi(y_j^{(i)};\X^{(i-1)})}$, where
		\[
			\chi(\cdot;\X^{(i-1)}) = \frac{1}{M}\sum_{j=1}^M \nu(\cdot;x_j^{(i-1)}).
		\]

		Resample: $(\W^{(i)},\Y^{(i)}) \rightarrow (\frac{1}{M}\mathbf{1}, \X^{(i)})$.\label{algline:ETAIS_resample}\;
	}
	Output $(\W, \Y)$.\;
\caption{The ensemble transform adaptive importance sampler (ETAIS).\label{alg:ETAIS}}
\end{algorithm}
\end{table}

Since the resampling does not give us a statistically identical sample
to that which is input, we cannot assume that the samples $\X^{(i)}$
are samples from the posterior. Therefore, as with serial
importance samplers, the weighted samples
$(\W,\Y)_{i=1}^N$ are the samples from the posterior that
we will analyse.

In each iteration of the algorithm, we are required to compute the
likelihood for each member of the ensemble. In the case where the
likelihood is expensive to compute, for example because it requires
the numerical approximation of an ODE or PDE, and/or there is a very
large amount of observations in the data, it will be expedient to
split this computational effort across multiple cores. Whether this is
efficient will depend on the architecture of the machine, and the
relative cost of communication across the cores, to the cost of the
likelihood evaluation. Evaluation of the denominator in the weights
also contributes to the overhead costs of this approach over and above
a standard MCMC method, but this too can be parallelised efficiently,
depending on communication costs of the architecture of the machine,
reducing the $\mathcal{O}(M^2)$ cost to $\mathcal{O}(M)$ over $M$
cores, for example.

The key is to choose a suitable transition kernel $\nu$ such that
if $\mathbf{X}^{(i)}$ is a representative sample of the posterior,
then the mixture density $\chi(\cdot;\X^{(i)})$ is a good
approximation of the posterior distribution. If this is the case,
the newly proposed states $\Y^{(i)}$ will also be a good sample of the
posterior with low variance in the weights $\W^{(i)}$.

In Section \ref{Sec:Num}, we will demonstrate how the algorithm
performs, primarily using random walk (RW) proposals. We do not claim that this choice is optimal, but
is simply chosen as an example to show that sharing information across
ensemble members can improve on the original MH algorithm and lead to
tolerances being achieved in fewer evaluations of $\G$. This is important since if
the inverse problem being tackled involves computing
the likelihood from a very large data set, or where the likelihood
requires the numerical solution of a differential equation, this could lead to a
large saving of computational cost. We have observed that
using more complex (and hence more expensive) kernels $\nu$, does not
significantly improve the speed of convergence of the algorithm for
posteriors that we have considered, although kernels such as those
used in MALA can be more stable for certain problems\cite{Paul}.

Care needs to be taken when choosing kernel(s) for the proposal distribution to
ensure that the overall mixture is absolutely continuous with respect to
the posterior distribution. For example, in many inverse problems we
may be looking to find the value of certain physical parameters which
may be strictly positive, or even bounded on an interval. In
this case, proposal kernels $\nu$ should be picked which have the same
support as the target. We will see an example of an unknown parameter
being bounded in Section \ref{sec:mixture}, where
$\beta$-distributed proposals are made to ensure that the distribution
is supported on $[0,1]$.

It is also possible, especially during the early stages of the
algorithm, where the mixture proposal is a poor approximation of the
target distribution, for a sample to be produced with a weight which
is orders of magnitude bigger than the rest. This is a problem, since
the resampling step will then lead to a vector of samples all
centred at the outlier sample. Fortunately these spikes in weights can
easily be detected, and in such a case, the sample for this iteration
can be removed. The ETPF can then be used to formulate a new mixture
sample with one sample placed at the outlier point, and the others
distributed according to the sample positions before the problematic
iteration. This ensures that the region where the problem occurred is
better represented by the mixture, and prevents further spikes in this region.

%%%%%%

\section{Automated tuning of algorithm parameters}\label{Sec:adapt}

Efficient selection of scaling parameters in MCMC algorithms is
critical to achieving optimal mixing rates and hence achieving fast
convergence to the target density. It is well known that a scaling
parameter which is either too large or too small results in a Markov
chain with high autocorrelation. One aspect worthy of consideration
with the ETAIS, is finding an appropriate proposal kernel $\nu$ such
that the mixture distribution $\chi$ is a close approximation to the
posterior density $\pi$.

Most MCMC proposals have parametric dependence which
allows the user to control their variance. For example, in the RW
proposal $y = x + \beta \eta$, the parameter $\beta$ is the standard deviation
of the proposal distribution. Therefore the proposal distributions can
be tuned such that they are slightly over-dispersed. This tuning
can take place during the burn-in phase of the algorithm. Algorithms
which use this method to find optimal proposal distributions are known
as adaptive MCMC algorithms, and have been shown to be convergent
provided that they satisfy certain
conditions\cite{roberts2007coupling,roberts2009examples}.
Alternatively, a costly trial and error scheme with short runs of the
MCMC algorithm can be used to find an acceptable value of $\beta$.

Algorithms which use mixture proposals, e.g. ETAIS, must tune the
variance of the individual kernels within the proposal mixture. This
adaptivity during early iterations has some added benefits over and above
finding an optimal parameter regime for the algorithm. If the initial
value of the proposal variance is chosen to be very large, then the early
mode-finding stages of the algorithm are expedited. Adaptively
reducing the proposal variances to an optimal value then allows us to
explore each region efficiently. Using an ensemble of chains allows
quick and effective assessment of the value of the optimal scaling
parameter.

In many MCMC algorithms such as the Random Walk Metropolis-Hastings
(RWMH) algorithm, the optimal scaling parameter can be found by
searching for the parameter value which gives an optimal acceptance
rate, e.g. for near Gaussian targets the optimal rates are 23.4\% for
RWMH and 57.4\% for MALA\cite{roberts2001optimal}. This method
is not applicable to ETAIS so we must use other statistics to optimise
the scaling parameter. Section~\ref{sec:statistics} gives some
possible methods for tuning $\beta$.

%%%%%%

\subsection{Statistics for Determining the Optimal Scaling Parameter}\label{sec:statistics}

\subsubsection{Determining optimal scaling parameter using error analysis}

When available, an analytic form for the target distribution allows us
to assess the convergence of sampling algorithms to the target
distribution. Common metrics for this task include the relative error
between the sample moments and the target's moments, or the
relative $L^2$ error between the sample histogram and the target
density, $\pi(x|D)$. The relative error in the $n$-th moment,
$\hat{m}_n$, is given by:
\begin{equation}\label{eq:34567}
	e = \left|\frac{\hat{m}_n - \mathbb{E}[X^n]}{\mathbb{E}[X^n]}\right|, \quad \text{where} \quad \hat{m}_n = \frac{1}{N}\sum_{i=1}^N \! x_i^n,
\end{equation}
and $\{x_i\}_{i=1}^N$ is a sample of size $N$.

The relative $L^2$ error, $E$, between a continuous density function to a
piecewise constant approximation of that density, can be defined by considering the
difference in mass between the self-normalised histogram of the
samples and the posterior distribution over a set of disjoint sets or
``bins'':
\begin{equation}\label{eqn:L2_error}
	E^2 = \sum\limits_{i=1}^{n_b}\left[\displaystyle\int_{R_i} \! \pi(s|D) \, \mbox{d}s - vB_i\right]^2 \Big/ \sum\limits_{i=1}^{n_b}\left[\displaystyle\int_{R_i} \! \pi(s|D) \, \mbox{d}s\right]^2,
\end{equation}
where the regions $\{R_i\}_{i=1}^{n_b}$ are the $d$-dimensional
histogram bins, so that $\bigcup_i R_i \subseteq X$ and
$R_i\cap R_j=\emptyset$, $n_b$ is the number of bins, $v$ is the
volume of each bin, and $B_i$ is the value of the $i$th bin. This
metric converges to the standard definition of the relative $L^2$
error as $v\rightarrow 0$.

These statistics cannot be used in general to find optimal values of
$\beta$ since they require the analytic solution, and are expensive to
approximate. However they can be used to assess the ability of other
indicators to find the optimal scaling parameters in a controlled
setting.

%%%%%%

\subsubsection{The effective sample size}\label{sec:ess}

The effective sample size, $\neff$, can be used to assess the
efficiency of importance samplers. Ideally, in each iteration, we
would like all $M$ of our samples to provide us with new information
about the posterior distribution. In practice, we cannot achieve a
perfect effective sample size of $M$.

The effective sample size of a weighted sample is defined in the
following way:
\[
	\neff = \frac{\left(\sum_{i=1}^M \! w_i\right)^2}{\sum_{i=1}^M \! w_i^2} \approx \frac{M\E(w)^2}{\E(w^2)} = M\left(1 - \frac{\mbox{var}(w)}{\mathbb{E}(w^2)}\right).
\]
The second two expressions are true when $M\rightarrow\infty$. From
the last expression, if the variance of the weights is zero then
$\neff = M$; this is our ideal scenario. In the limit
$M\rightarrow\infty$, maximising the effective
sample size is equivalent to minimising the variance of the weights.

The statistic $\neff$ is easier to deal with than the variance of the
weights, as it takes values on $[1, M]$ while the variance of the
weights takes values on $\mathbb{R}^+$. Moreover the variance of
the weights can vary
over many orders of magnitude causing numerical instabilities, so that the effective
sample size is more desirable as an indicator of optimality.

In this paper we tune our simulations using the effective sample size
statistic. We calculate this statistic using a sample size of $Mn_k$,
where $1 \leq n_k\leq N$ is sufficiently large enough to obtain a
reasonable estimate of $\neff$ with scaling parameter $\delta_k$.
Calculating $\neff$ over these subsets of the simulations tends to
underestimate the optimal value of the scaling parameter due to the
possibility of missing unlikely proposals with extreme weights.

The effective sample size also has another useful property; if we
consider the algorithm in the early stages, for example, we have all
$M$ samples in the tails of the target searching for a mode. The
sample closest to the mode will have an exponentially higher weight
and the effective sample size will be close to 1. In later stages the
ensemble populates high density regions, and better represents the
posterior distribution. This leads to smaller variation in the
weights, and a higher effective sample size. By looking for
approximations of the effective sample size which look like a
stationary distribution, we can tell when the algorithm is working
efficiently.

\begin{figure}[htb]
\centering
\includegraphics[width=0.65\textwidth]{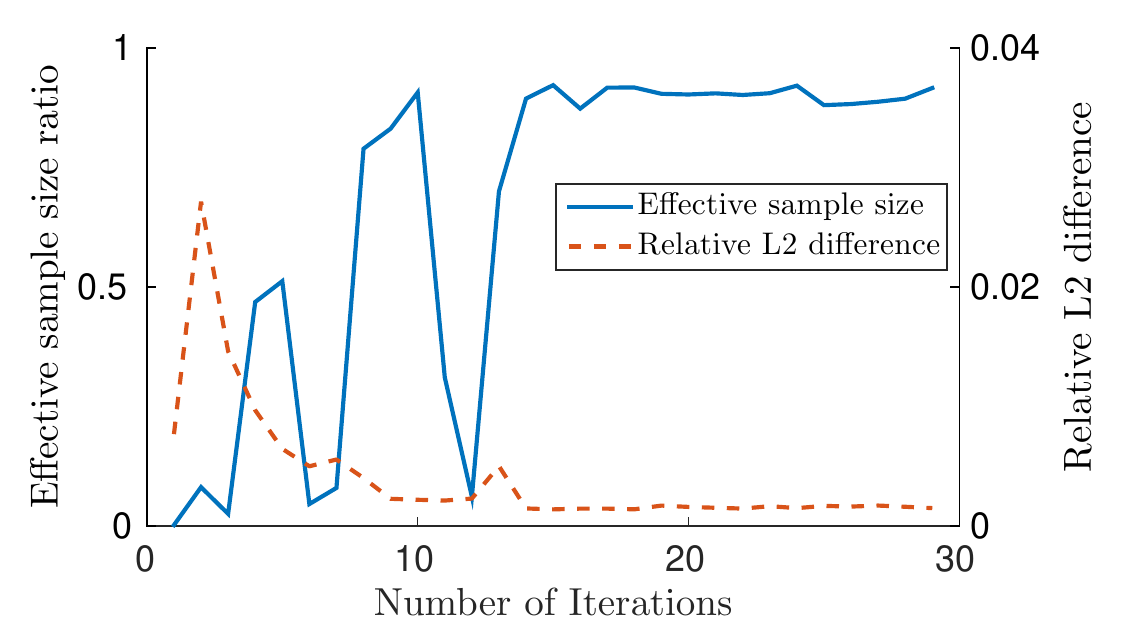}
\caption{The effective sample size ratio and relative $L^2$ difference
$E$ during the first 30 iterations. These numerics are taken from the
example in Section~\ref{sec:mixture_conv} using the ETAIS-RW algorithm.}
\label{fig:neff-burnin}
\end{figure}

Figure~\ref{fig:neff-burnin} demonstrates that the effective sample
size flattens out as the relative $L^2$ difference between the
posterior distribution and the proposal distribution stabilises close
to its minimum.

\begin{figure}[htb]
\centering
\subfigure[Contours showing optimal ranges of the scaling parameter.]{\includegraphics[width=0.47\textwidth]{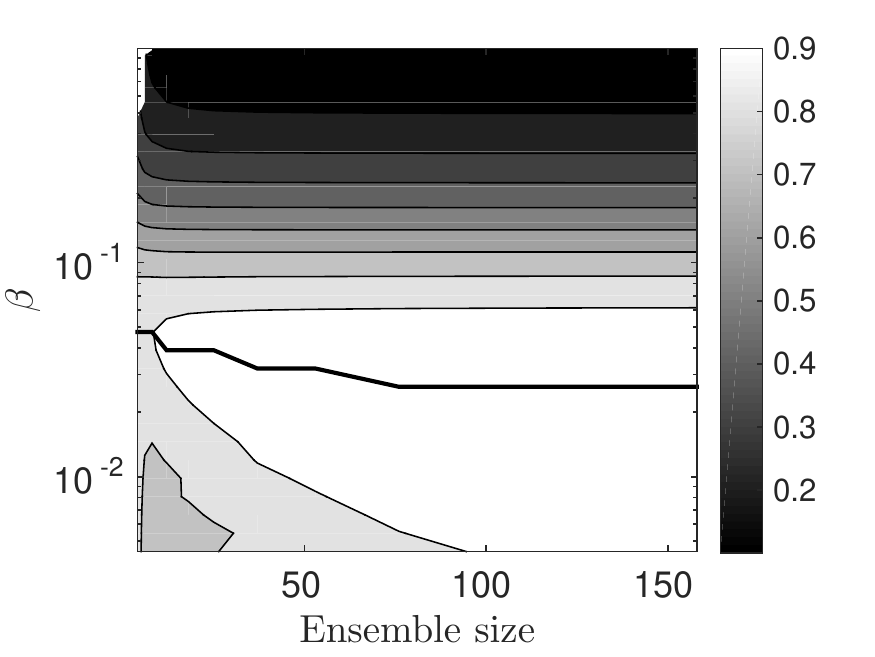}}
\subfigure[The highest effective sample size achieved for each ensemble size.]{\includegraphics[width=0.47\textwidth]{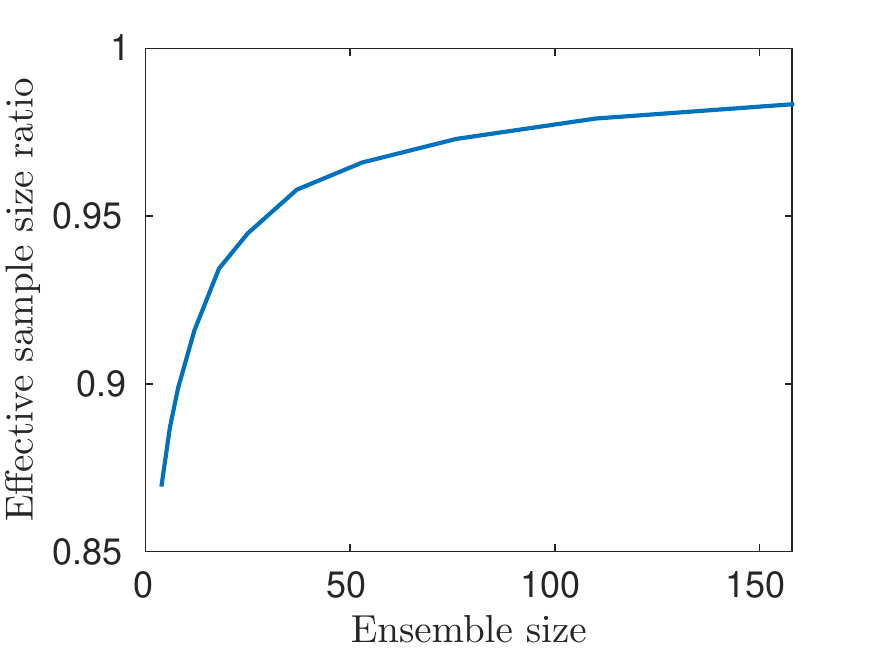}}
\caption{The behaviour of the effective sample size as the ensemble
  size increases, considering the example in Section~\ref{sec:problem 1} using the ETAIS-RW algorithm.}
\label{fig:neff-M}
\end{figure}

Figure~\ref{fig:neff-M} shows how the effective sample size used in
ETAIS is affected by the ensemble size. We see from subfigure (a) that
as the ensemble size increases, the optimal scaling parameter
decreases. This is expected since the larger ensemble allows for
finer resolution in the proposal distribution approximation of the
posterior distribution. We also
see that the algorithm becomes less sensitive to changes in the
scaling parameter as the ensemble size increases. Subfigure (b) shows
that as the ensemble size increases, the algorithm becomes more
efficient.

\subsection{Adaptively Tuned ETAIS}\label{sec:adapt}

To adapt the scaling parameter $\beta$, we use a version of the
gradient ascent method modified for a stochastic function. Some more
sophisticated examples are described in
\cite{roberts2009examples,Ji2013adaptive,andrieu2006ergodicity}.

\begin{table}[!ht]
\centering
\begin{algorithm}[H]
\DontPrintSemicolon
\BlankLine
	Define update times $\{n_k\}_k$.\;
	\For{$n=1,\dots,N$}{
		Complete steps \ref{algline:ETAIS_propose}-\ref{algline:ETAIS_resample} of Algorithm~\ref{alg:ETAIS}.\;
		\If {$n \in \{n_k\}$}
		{
			Divide ensemble into two halves. Use these halves to estimate the gradient in $\neff$ at $\beta_k$.\;
			Update $\beta_k$ using gradient ascent,
				\[
					\beta_{k+1} = \beta_k + \gamma\nabla\neff.
				\]
			\label{algline:gradient_ascent}
		}
	}
\caption{Adaptively tuned ETAIS algorithm.\label{alg:adaptETAIS}}
\end{algorithm}
\end{table}
Our adaptive algorithm is given in Algorithm~\ref{alg:adaptETAIS}. We choose update times which, as suggested in Section~\ref{sec:ess} allow for a reasonable estimate of the effective sample size, but do not waste too many iterations. In Step~\ref{algline:gradient_ascent}, the gradient ascent parameter $\gamma$ may decrease over time, e.g. as a function of $n_k$.

\begin{figure}[!htb]
\centering
\subfigure[AETAIS-RW algorithm adapting variance.]{\includegraphics[width=0.495\textwidth]{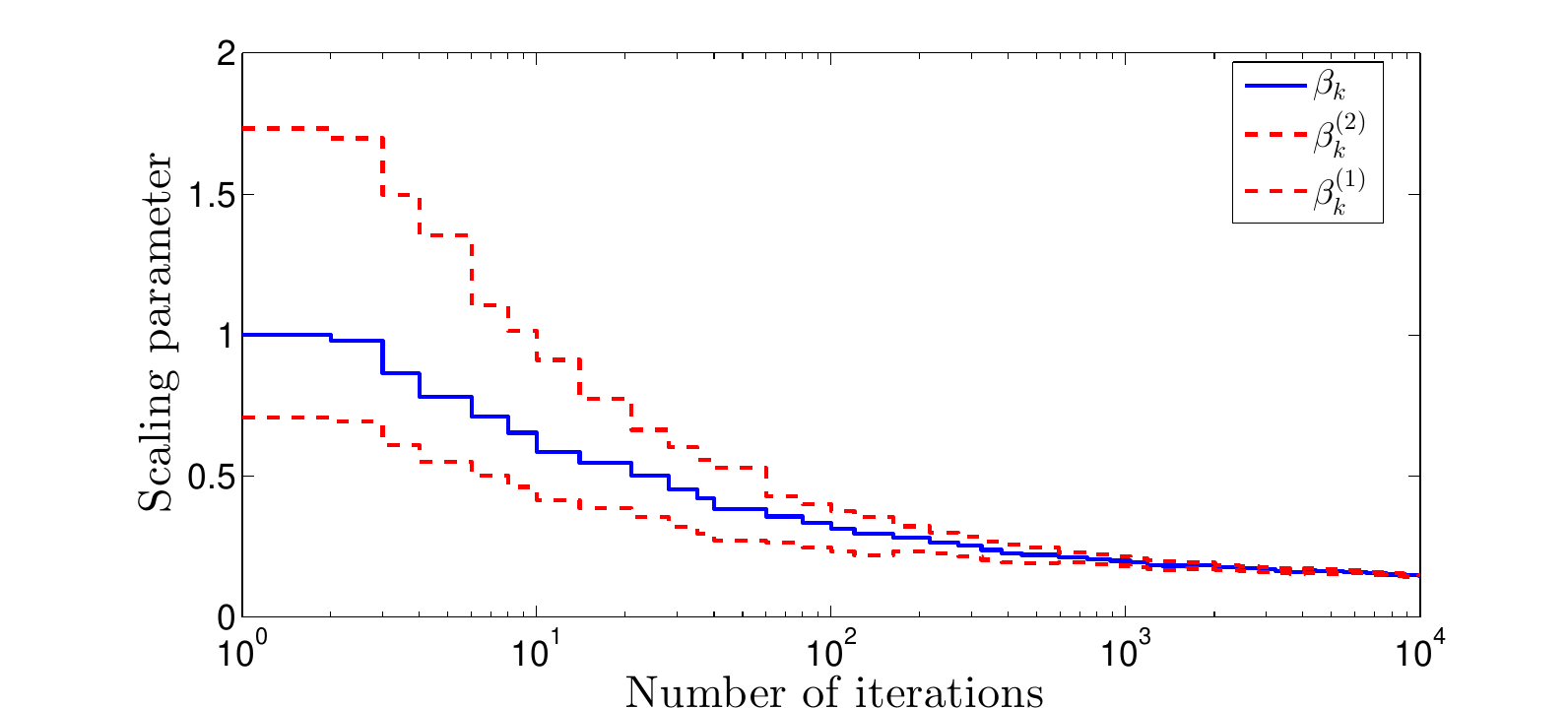}}
\subfigure[AETAIS-RW algorithm maximising ESS.]{\includegraphics[width=0.495\textwidth]{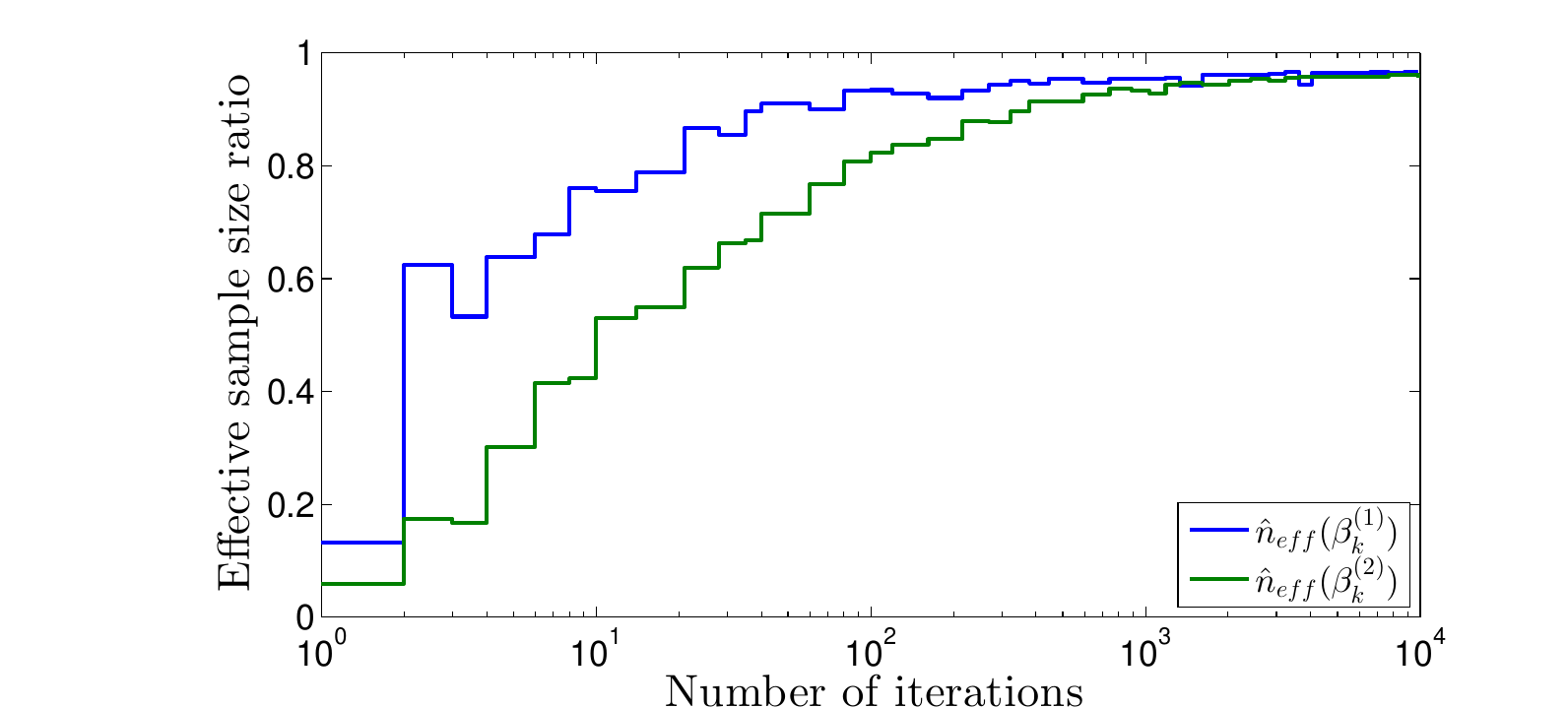}}
\subfigure[AETAIS-MALA algorithm adapting variance.]{\includegraphics[width=0.495\textwidth]{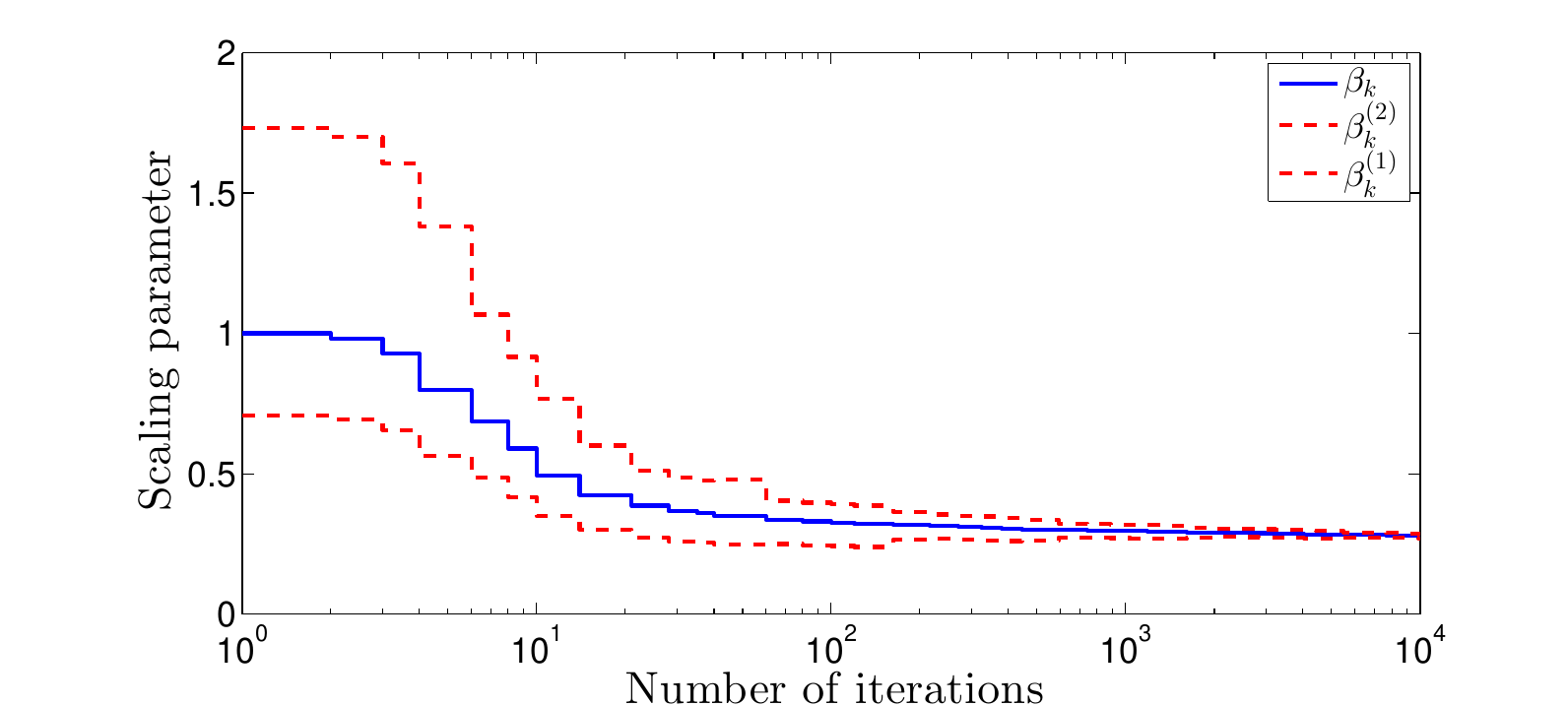}}
\subfigure[AETAIS-MALA algorithm maximising ESS.]{\includegraphics[width=0.495\textwidth]{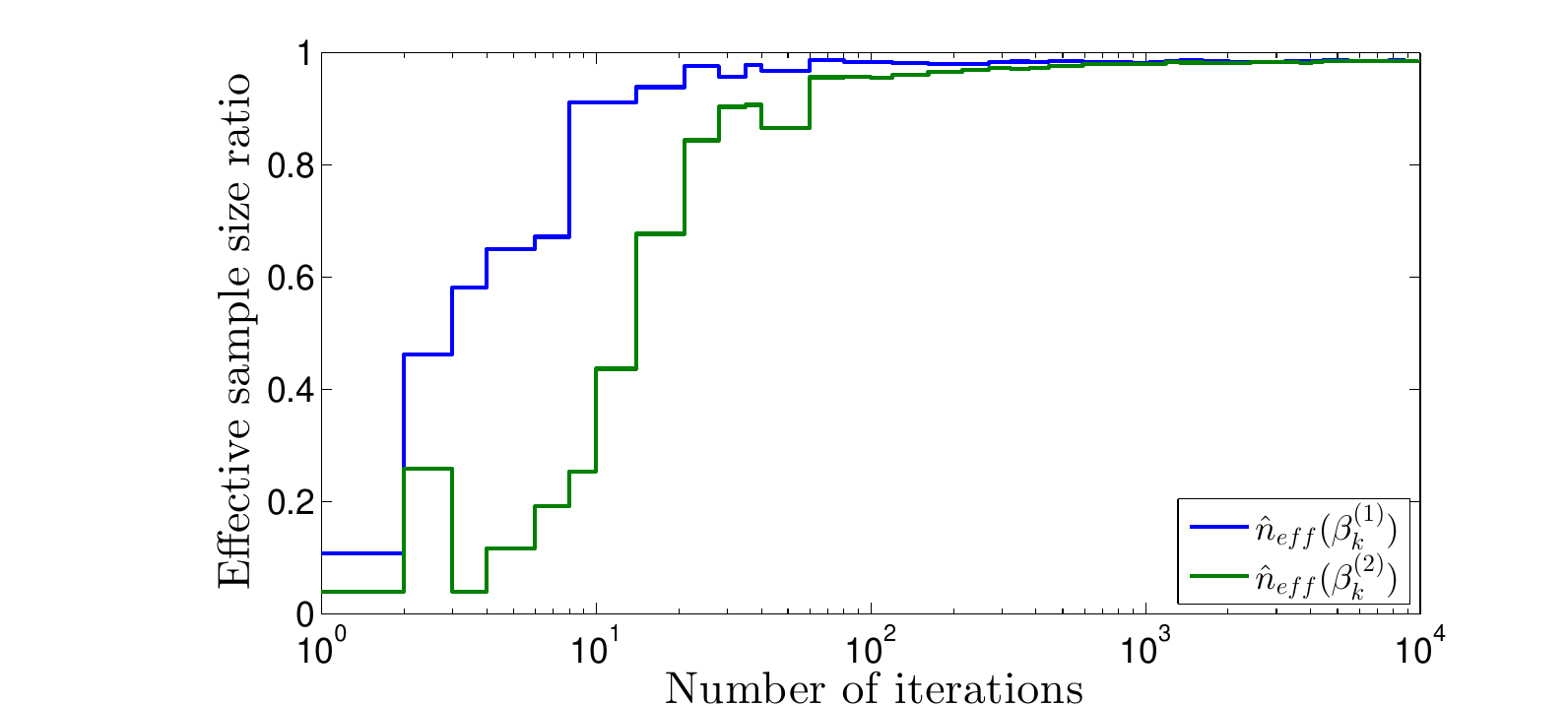}}
\caption{Demonstration of the convergence of the adaptive scaling parameter in the adaptively tuned ETAIS algorithms.}
\label{fig:G1_adapt_trace_ETAIS}
\end{figure}

In the context of an expensive likelihood evaluation, slow convergence
of the algorithmic parameters in this adaptive regime would be a
concern. However, this is rarely an issue, due in part to the speed at
which the ETAIS algorithm completes the burn-in phase. Figure
\ref{fig:G1_adapt_trace_ETAIS} demonstrates this for two different
versions of the ETAIS methodology with a RW proposal and a MALA\cite{roberts1998optimal}
proposal. This picture is typical in our experience, with the scaling
parameters converging in $\mathcal{O}(10^2)$ iterations to an optimal regime.
%%%%%%

\section{Multinomial Transformation}\label{sec:MT}

Although the ETPF produces the optimal linear coupling, it can also become quite costly as the number of ensemble members
is increased. It is arguable that in the context of ETAIS, we do not
require this degree of accuracy, and that a faster more approximate
method for resampling could be employed. One approach would be to use
the bootstrap resampler, which simply takes the $M$ ensemble members'
weights and constructs a multinomial distribution, from which $M$
samples are drawn. This is essentially the cheapest resampling
algorithm that one could construct. However it too has some
drawbacks. The algorithm is random, and as such it is possible for all
of the ensemble members in a particular region not to be sampled. This
could be particularly problematic when attempting to sample from a
multimodal distribution, where it might take a long time to find one
of the modes again. The bootstrap filter is also not guaranteed to
preserve the mean of the weighted sample, unlike the ETPF.

Ideally, we would like to use a resampling algorithm which is not
prohibitively costly for moderately or large sized ensembles,
which preserves the mean of the samples, and which makes it much
harder for the new samples to forget a significant region in the
density. This motivates the following algorithm, which we refer to as the
multinomial transformation (MT), which is a greedy
approximation of the ETPF resampler.

Instead of sampling $M$ times from an $M$-dimensional multinomial
distribution as is the case with the bootstrap algorithm, we sample
once each from $M$ different multinomials. Suppose that we have $M$
samples $y_n$ with weights $w_n$. The multinomial sampled from in the
bootstrap filter has a vector of probabilities given by:
\begin{equation*}
\frac{1}{\sum w_n} [w_1,w_2,\ldots,w_M] = \bar{\bf w},
\end{equation*}
with associated states $y_n$.
We wish to find $M$ vectors $\{{\bf p}_1,{\bf p}_2,\ldots,{\bf p}_M\}
\subset \mathbb{R}^M_{\ge 0}$
such that  $\frac{1}{M} \sum {\bf p}_i = \bar{\bf w}$. The MT is then
given by a sample from each of the multinomials defined by the vectors
${\bf p}_i = [p_{i,1},p_{i,2},\ldots,p_{i,M}]$ with associated states ${\bf y}_i$. Alternatively, as with the ETPF, a deterministic sample
can be chosen by picking each sample to be equal to the mean value of
each of these multinomial distributions, i.e. each new sample
$\hat{x}_i$ is given by:
\begin{equation}
\hat{x}_i = \sum p_{i,j} x_j, \qquad i \in \{1,2,\ldots,M\}.
\end{equation}

The resulting sample has several properties which are advantageous in
the context of being used with the ETAIS algorithm. Firstly, we have
effectively chopped up the multinomial distribution used in the
bootstrap filter into $M$ pieces, and we can guarantee that exactly
one sample will be taken from each section. This leads to a much
smaller chance of losing entire modes in the density, if each of the
sub-multinomials is picked in an appropriate fashion. Secondly, if we do not make a random sample for
each multinomial with probability vector $\bf p_i$ but instead take
the mean of the multinomial to be the sample, this algorithm preserves
the mean of the sample exactly. Lastly, as we will see shortly, this
algorithm is significantly less computationally intensive than the
ETPF.

There are of course infinitely many different ways that one could use
to split the original multinomial up into $M$ parts, some of which
will be far from optimal. The method that we have chosen is loosely
based on the idea of optimal transport. We search out states with the
largest weights, and choose a cluster around these points based on
the closest states geographically. This method is not optimal since
once most of the clusters have been selected the remaining states
may be spread across the parameter space.

\begin{table}[!ht]
\centering
\begin{algorithm}[H]
\DontPrintSemicolon
\BlankLine
	$\b{z} = M\bar{\b{w}}$.\;
	\For {$i = 1,\dots, M$}
	{
		$J = \argmax_j z_j$.\;
		$p_{i,J} = \min\{1,z_J\}$.\;
		$z_J = z_J - p_{i,J}$.\;
		\While {$\sum_j p_{i,j} <1$}
		{
			$K = \argmin_{k \in \{k|z_k>0\}} \|y_J - y_k\|$.\;
			$p_{i,K} = \min\{1-\sum_j p_{i,j}, z_K\}$.\;
			$z_K = z_K - p_{i,K}$.\;
		}
		$x_i = \sum_k p_{i,k}y_k$.\;
	}
\caption{The  multinomial transformation (MT).\label{alg:MT}}
\end{algorithm}
\end{table}

Algorithm~\ref{alg:MT} describes the basis of the algorithm with
deterministic resampling, using the means of each of the
sub-multinomials as the new samples. This resampler was designed with the aims of being
numerically cheaper than the ETPF, and more accurate than straight multinomial
resampling. Therefore we now present numerical examples which
demonstrate this.

\begin{figure}[htb]
\centering
\subfigure[Relative error in $\mathbb{E}(X)$]{\includegraphics[width=0.45\textwidth]{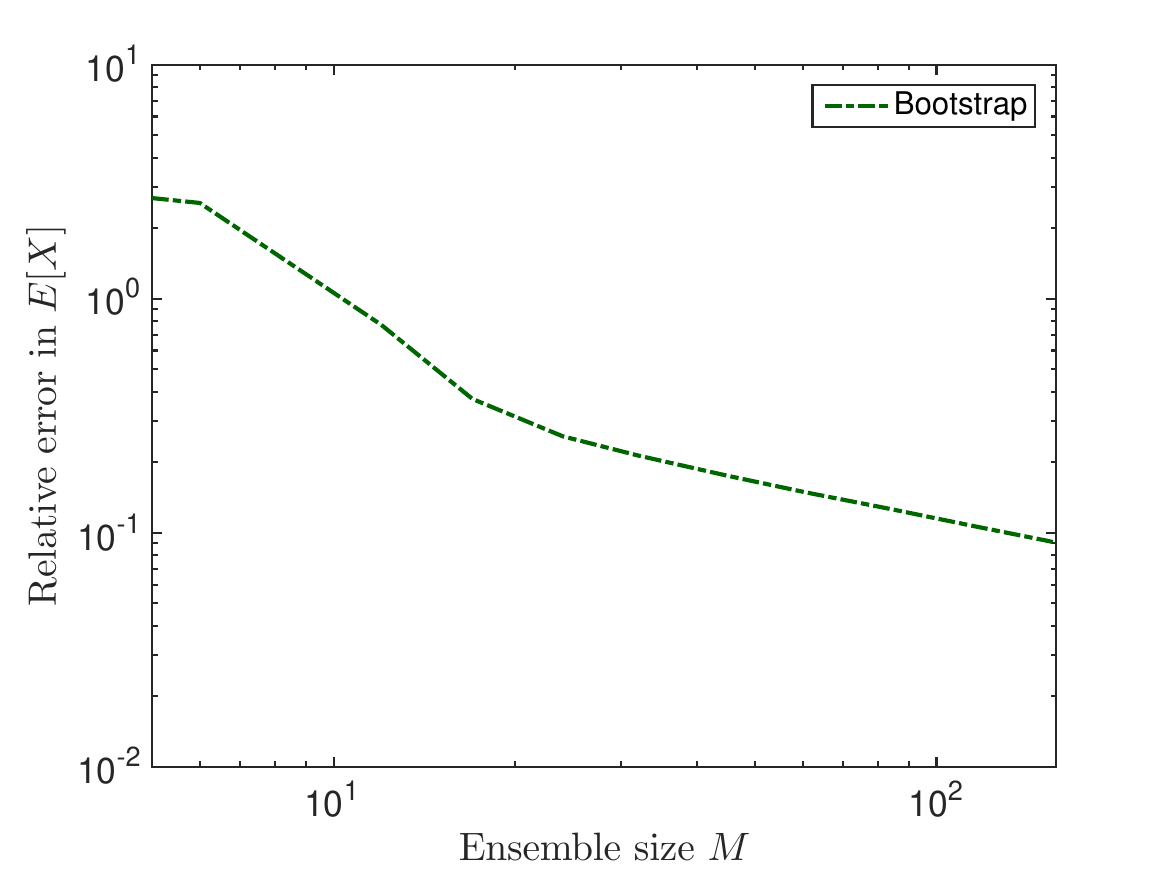}}
\subfigure[Relative error in
$\mathbb{E}(X^2)$]{\includegraphics[width=0.45\textwidth]{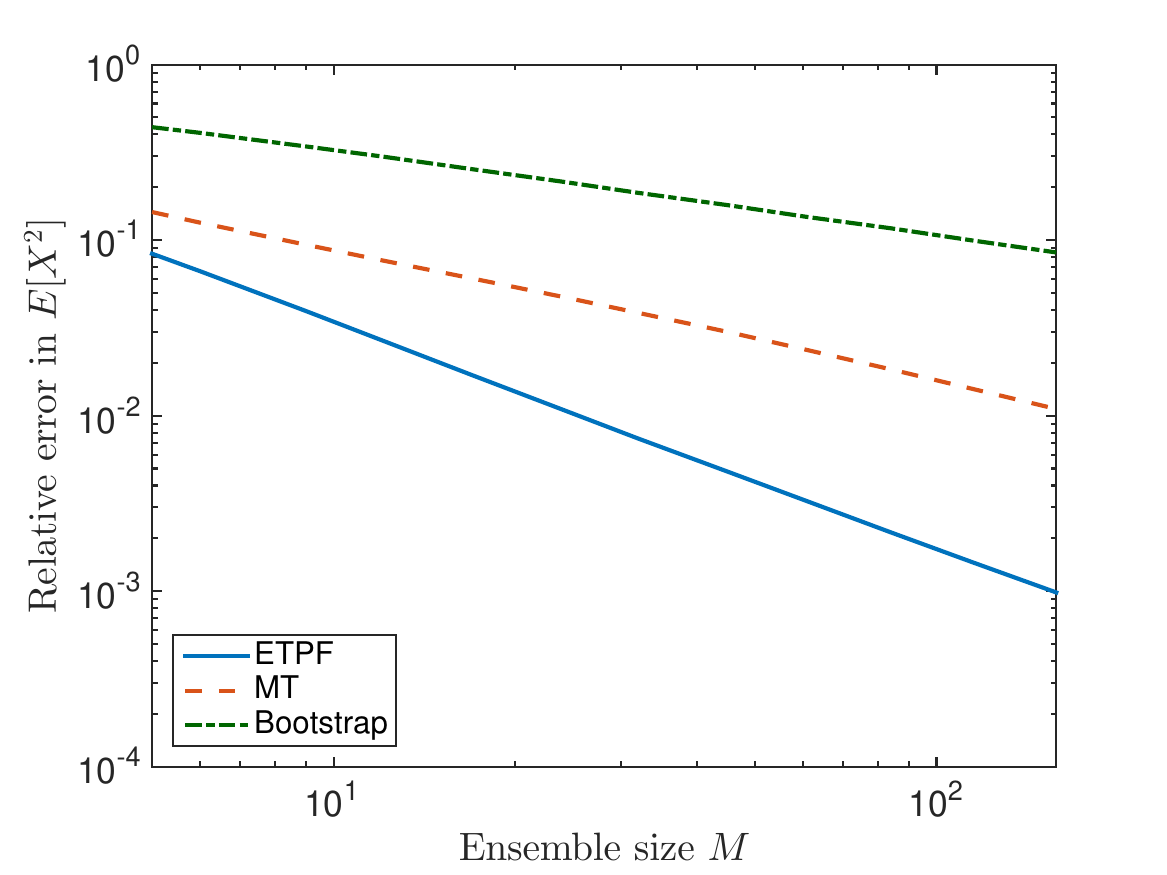}}\\
\subfigure[Relative error in $\mathbb{E}(X^3)$]{\includegraphics[width=0.45\textwidth]{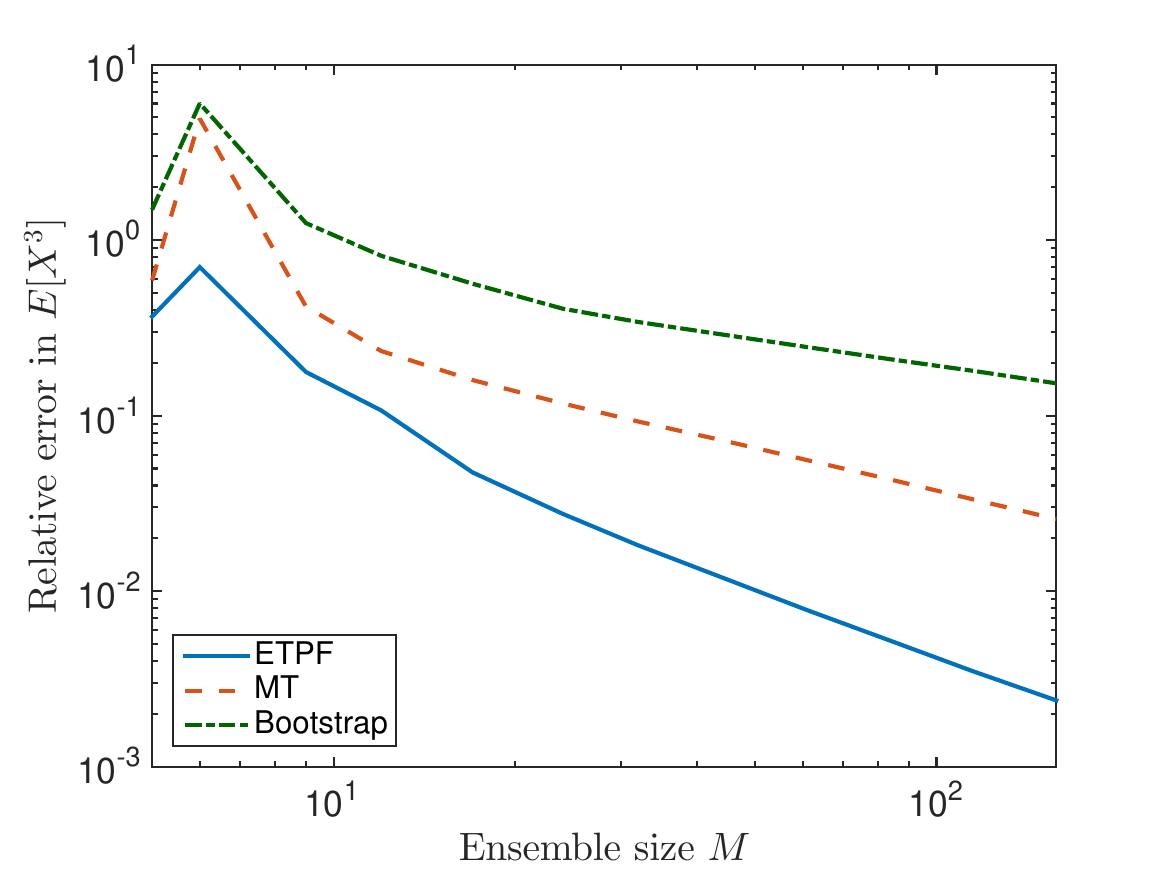}}
\subfigure[Cost of resamplers per iteration]{\includegraphics[width=0.45\textwidth]{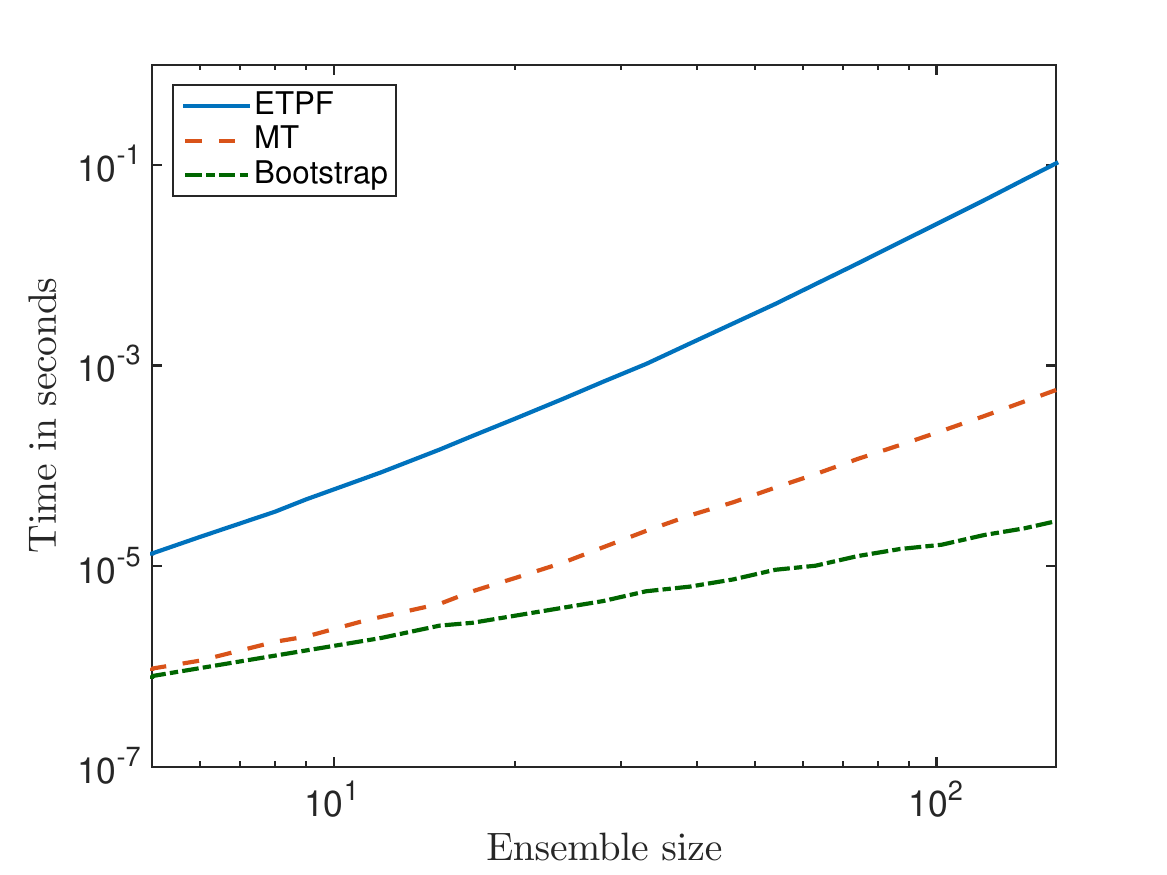}}\\
\caption{Comparison of the performance between different resampling schemes. The example in Section~\ref{sec:problem 1} is implemented for this demonstration.}
\label{fig:MT}
\end{figure}

To test the accuracy and speed of the three resamplers (ETPF,
bootstrap and MT), we drew a sample of size $M$ from the proposal distribution
$\mathcal{N}(1,2)$. Importance weights were assigned, based on a
target distribution of $\mathcal{N}(2,3)$. The statistics of the
resampled outputs were compared with the original weighted samples. Figure \ref{fig:MT} (a)-(c) show how the relative errors in the first
three moments of the samples changes with ensemble size $M$ for the three different
resamplers. As expected, the MT lies somewhere between the high
accuracy of the ETPF and the less accurate bootstrap
resampling. Note that only the error for the bootstrap multinomial
sampler is presented for the first moment since both the ETPF and the
MT preserve the mean of the original weighted samples up to machine precision. Figure \ref{fig:MT} (d) shows how the computational cost,
measured in seconds, scales with the ensemble size for the three
different methods, where timings have been taken from simulations on a Dell
server with four 8-core 3.3GHz CPUs and 64Gb memory. These results demonstrate that the MT behaves how we wish, and
importantly ensures that exactly one sample of the output will lie in
each region with weights up to $\frac{1}{M}$ of the total.

We will use the MT in the numerics in Section \ref{sec:mixture}
where we have chosen to use a larger ensemble size. We do not claim
that the MT is the optimal choice within ETAIS, but it does have
favourable features, and demonstrates how different choices of
resampler can affect the speed and efficiency of the ETAIS algorithm.

%%%%%%

\section{Consistency of ETAIS}\label{sec:consistency}
As with all importance sampling schemes, we must have absolute
  continuity of target distribution with respect to the proposal
  distribution if we wish to achieve convergence. This can usually be
  ensured by picking the proposal kernels $\nu$ to be from the same
  distribution type as the prior distribution.
  
As outlined in~\cite{martino2015adaptive} consistency of population AIS algorithms can be considered in two different cases. In the first case we fix the number of iterations $N < \infty$, but allow the population size to grow to infinity $M\rightarrow\infty$. In the second case we hold the population size fixed $M<\infty$, and allow infinitely many iterations $N\rightarrow\infty$.

In case one, from standard importance sampling results we know that for an iteration $n$, as $M\rightarrow\infty$, we obtain a consistent estimator for any statistic of interest,
\[
	\hat{\phi}(\Theta) \approx \frac{1}{\hat{Z}}\sum\limits_{i=1}^M \! \frac{1}{M}w_i\phi(\theta_i) \rightarrow \phi(\Theta),
\]
where the normalisation constant, $\hat{Z} = \frac{1}{M}\sum_{i=1}^M \! w_i$, also converges to the true normalisation constant $Z$~\cite{robert2013monte}.

Case two is slightly more involved. Estimation of the normalisation constant $Z$ is biased, and so estimates of statistics are sums of independent but biased estimators. Since the estimators are independent, the proof of consistency of ETAIS in this second limit can be approached in the same way as the pMC algorithm, where it has been shown that $\hat{Z}\rightarrow Z$~\cite{robert2013monte}. Since the normalisation constant is consistent, sums of the independent estimators are also consistent.

%%%%%%

\section{A useful property of the ETAIS algorithm for multimodal
  distributions}\label{sec:useful}

\begin{figure}[!ht]
\begin{center}
\includegraphics[width=\textwidth]{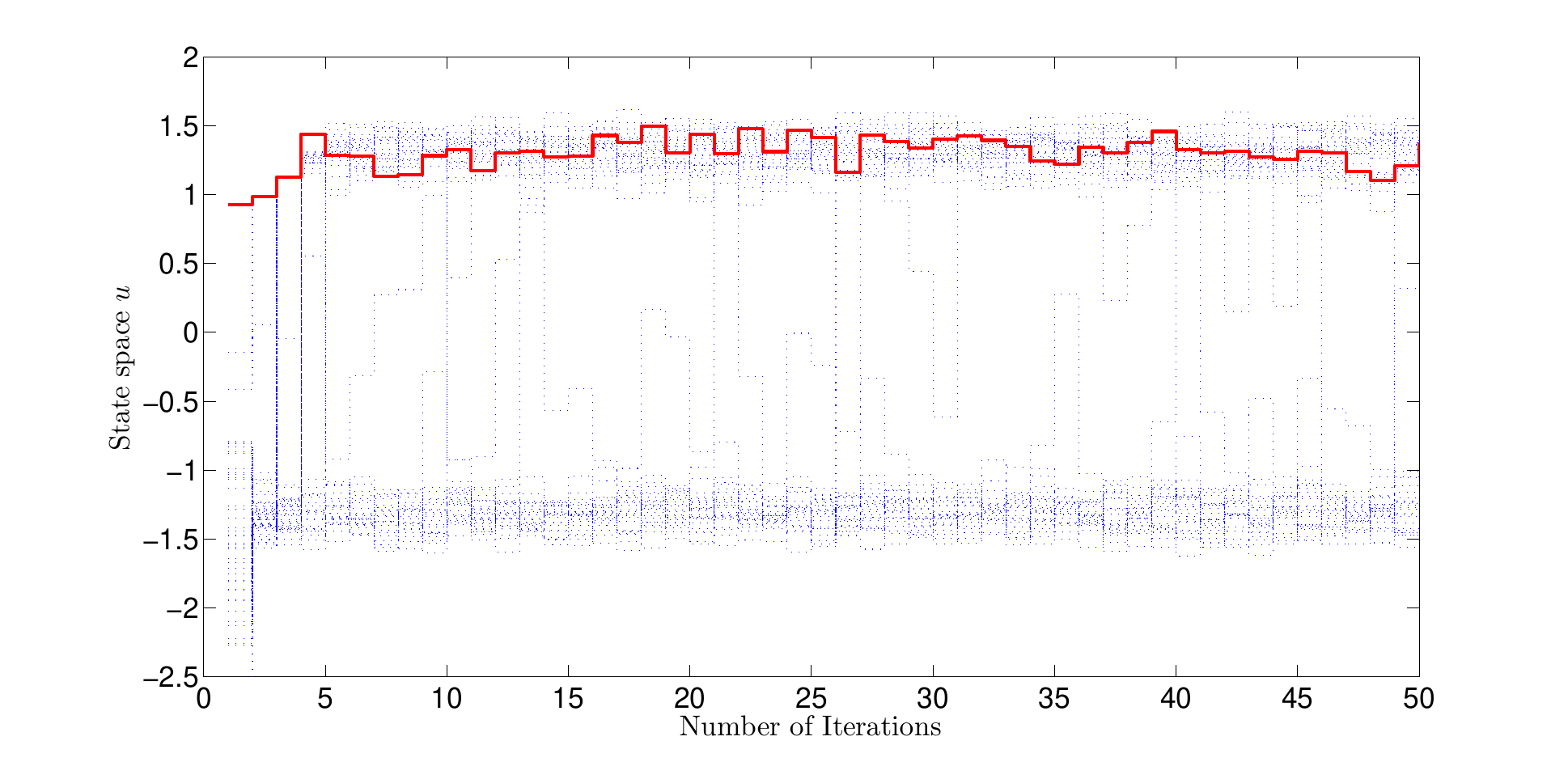}
\caption{This figure demonstrates the redistribution property of the
ETAIS algorithm for a bimodal target with equally weighted modes. Initially there is one chain in the positive mode, and
49 chains in the negative mode.}
\label{fig:BM2_suction}
\end{center}
\end{figure}

The biggest issue for the Metropolis-Hastings algorithms when sampling
from multimodal posterior distributions is that frequent switches
between the modes are required if the method is to converge
efficiently. The ETAIS algorithm tackles this problem with its
resampling step. The algorithm uses its dynamic mixture proposal distribution to build up an
approximation of the posterior at each iteration, and then compares
this to the posterior distribution via the weights function. Any large
discrepancy in the approximation will result in a large or small
weight being assigned to the relevant chain, meaning the chain will
either pull other chains towards it or be sucked towards a chain with
a larger weight. In this way, the algorithm allows chains to
`teleport' to regions of the posterior which are in need of more
exploration. Figure~\ref{fig:BM2_suction} shows the trace of a
simulation of the ETAIS-RW algorithm for a bimodal example with
equally weighted modes, with initially 1 chain in the
positive mode, and 49 chains in the negative mode. It takes only a
handful of iterations for the algorithm to balance out the chains into
25 chains in each mode. The chains switch modes without having to
climb the energy gradient in the middle.

\section{Numerical Examples}\label{Sec:Num}
The numerical examples in this section were all computed on a single
core of a Dell
server with four 8-core 3.3GHz CPUs and 64Gb memory. Throughout we
measure the efficiency of the algorithms in terms of the number of
iterations (i.e. likelihood evaluations) which are required in order
to reach a given order of accuracy. The reason for this is that this
method is designed with a particular type of challenging inverse
problem in mind, namely one which is low dimensional, has an expensive
likelihood (such as a PDE solve, and/or very large data) which dwarfs
the overhead cost of ETAIS, and which could have complex posterior
structure, such as correlations, multimodality or other highly
non-Gaussian features.

The numerics contained herein compare the ETAIS approach that we
  have presented against Metropolis-Hastings methods with
  the same type of proposal as is used by the kernels within
  ETAIS. Since we propose that ETAIS is a potential framework for
  parallelised Bayesian computations, we compare ETAIS with $N$
  ensemble members against an independent ensemble of $N$
  Metropolis-Hastings chains. All implementations were
  computed in serial, but the results demonstrate a clear speed-up
  which could be further exploited through parallelisation.

% In this section we demonstrate convergence of the ETAIS algorithm. We
% also investigate the convergence properties of ETAIS compared with
% naively parallelised Metropolis-Hastings samplers. We assess the
% performance of the ETAIS algorithm using the relative $L^2$ error
% defined in Equation~\eqref{eqn:L2_error}, as well as the relative error in the
% first moment (Equation~\eqref{eq:34567}).
%%%%%%

% \subsection{Numerical implementation}
% \label{sec:Implementation P1}

% For each example, we perform the following three tasks.

% \begin{enumerate}
% 	\item {\bf Finding the optimal scaling parameters}: We chose 32 values of $\beta$ evenly spaced on a log scale in the interval $[10^{-5}, 2]$. We ran the ETAIS and MH algorithms for one million iterations, each with an ensemble size of $M=50$. We took 32 repeats and used the geometric means of the sample statistics in Section~\ref{sec:statistics} to find the optimal scaling parameter.
% 	\item {\bf Measuring convergence of nonadaptive algorithms}: We ran the ETAIS and MH algorithms, using the scaling parameters found in Step 1, for one million iterations, again with $M=50$. The simulations are repeated 32 times. The performance is evaluated using the relative $L^2$ error defined in Equation~\eqref{eqn:L2_error}.
% 	\item {\bf Measuring convergence of adaptive algorithms}: We run the
% adaptive algorithms under the same conditions as the nonadaptive
% algorithms, and again use the relative $L^2$ error to compare
% efficiency. The adaptive algorithms are initialised with $\beta_1$ = 1.
% \end{enumerate}

% %%%%%%

\subsection{Automated variance tuning}
\label{sec:problem 1}
In this first example, we target a simple 1D Gaussian distribution, in
order to show the approximate equivalence of optimising the variance
of the proposal kernels using the effective sample size, and the $L^2$
error, which of course is not available to us in practical
applications.

\begin{figure}[htb]
\centering
\subfigure[RWMH algorithm.]{\includegraphics[width=0.45\textwidth]{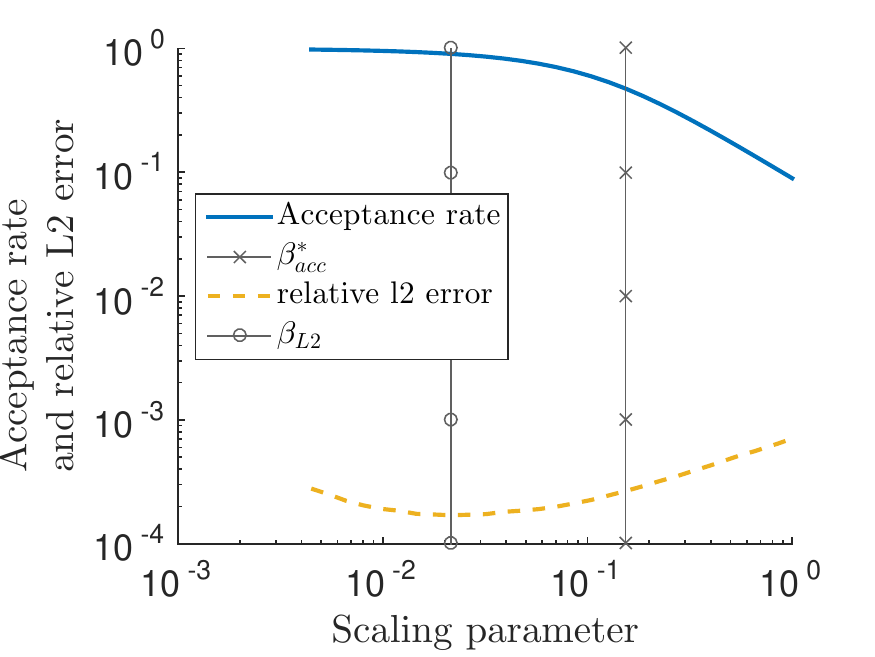}}
\subfigure[ETAIS-RW
algorithm.]{\includegraphics[width=0.45\textwidth]{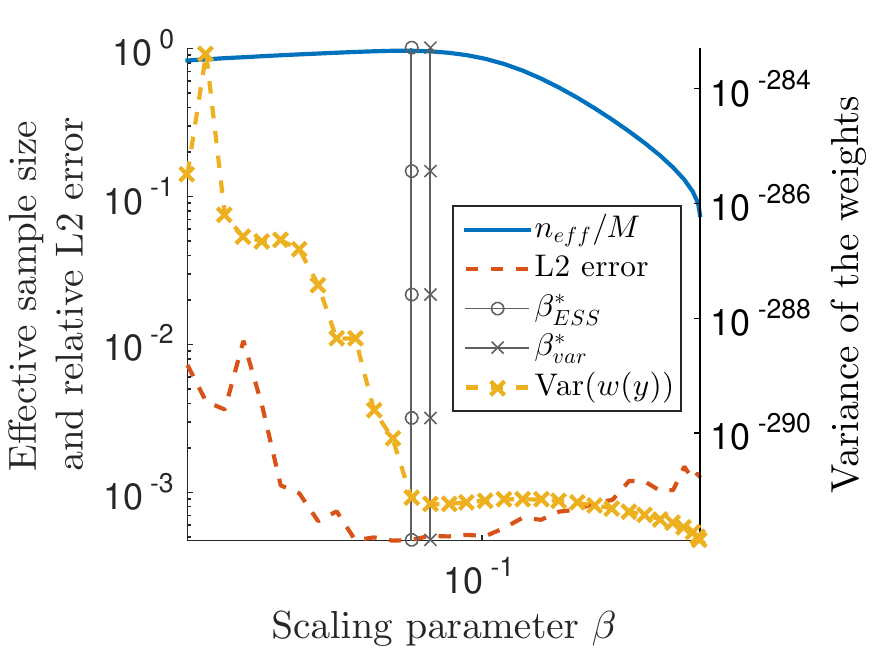}}

\caption{Finding optimal values of $\beta$ for a 1D Gaussian density. Bottom left: Convergence of the (A)RWMH and (A)ETAIS-RW algorithms.}
\label{fig:P1 opt beta}
\end{figure}

\begin{table}[!htb]
    \begin{minipage}{.5\linewidth}
      \centering
        \begin{tabular}{|c|r|}
	\hline
	Statistic											& RWMH \\ \hline
	$\beta_{\text{L2}}^*$								& 2.1e-2 \\
	$\beta_{\%}^*$									& 1.5e-1 \\
	Acceptance Rate ($\beta_{\text{L2}}^*$)				& 9.0e-1 \\
	Acceptance Rate ($\beta_{\%}^*$)					& 5.0e-1 \\
	\hline
	\end{tabular}
    \end{minipage}%
    \begin{minipage}{.5\linewidth}
      \centering
        \begin{tabular}{|l|r|r|}
	\hline
	Statistic							& ETAIS-RW \\ \hline
	$\beta_{\text{eff}}^*$				& 4.7e-2 \\
	$\beta_{\text{var}(w(y))}^*$		& 5.8e-2 \\
	$\beta_{\text{L2}}^*$				& 3.9e-2 \\
	\hline
	\end{tabular}
    \end{minipage}
	\vspace{1mm}
	\caption{Optimal values of $\beta$ summarised from Figure~\ref{fig:P1 opt beta}. Statistics calculated as described in Section~\ref{sec:statistics}. The values $\beta^*_{\text{L2}}$ and $\beta^*_{\%}$ are the optimal scaling parameters found by optimising the relative $L^2$ errors and acceptance rate respectively. Similarly $\beta_{\text{eff}}^*$ and $\beta_{\text{var}(w(y))}^*$ optimise the effective sample size and variance of the weights statistics.}
	\label{table:prob1 opt beta}
\end{table}

Figure~\ref{fig:P1 opt beta} (a) shows the two values of $\beta$ which
are optimal according to the acceptance rate and relative $L^2$ error
criteria for the RWMH algorithm. The smaller estimate comes from the
relative $L^2$ error, and the larger from the acceptance rate. The
results in Figure~\ref{fig:P1 opt beta} are summarised in
Table~\ref{table:prob1 opt beta}. Since in general we cannot calculate
the relative $L^2$ error, we must optimise the algorithm using the
acceptance rate. From the relative $L^2$ error curve we can see that
the minimum is very wide and despite the optimal values being very
different there is not a large difference in the convergence rate.

Figure~\ref{fig:P1 opt beta} (b) shows the effective sample size ratio
compared to the error analysis and the variance of the weights. The
relative $L^2$ error graph is noisy, but it is clear that the maximum
in the effective sample size and the minimum in the variance of the
weights are both close to the minimum in the relative $L^2$ error. Due
to this we say that the estimate of the effective sample size found by
averaging the statistic over each iteration is a good indicator for
the optimal scaling parameter.

\begin{figure}[htb]
\centering 
\subfigure[Convergence]{\includegraphics[width=0.45\textwidth]{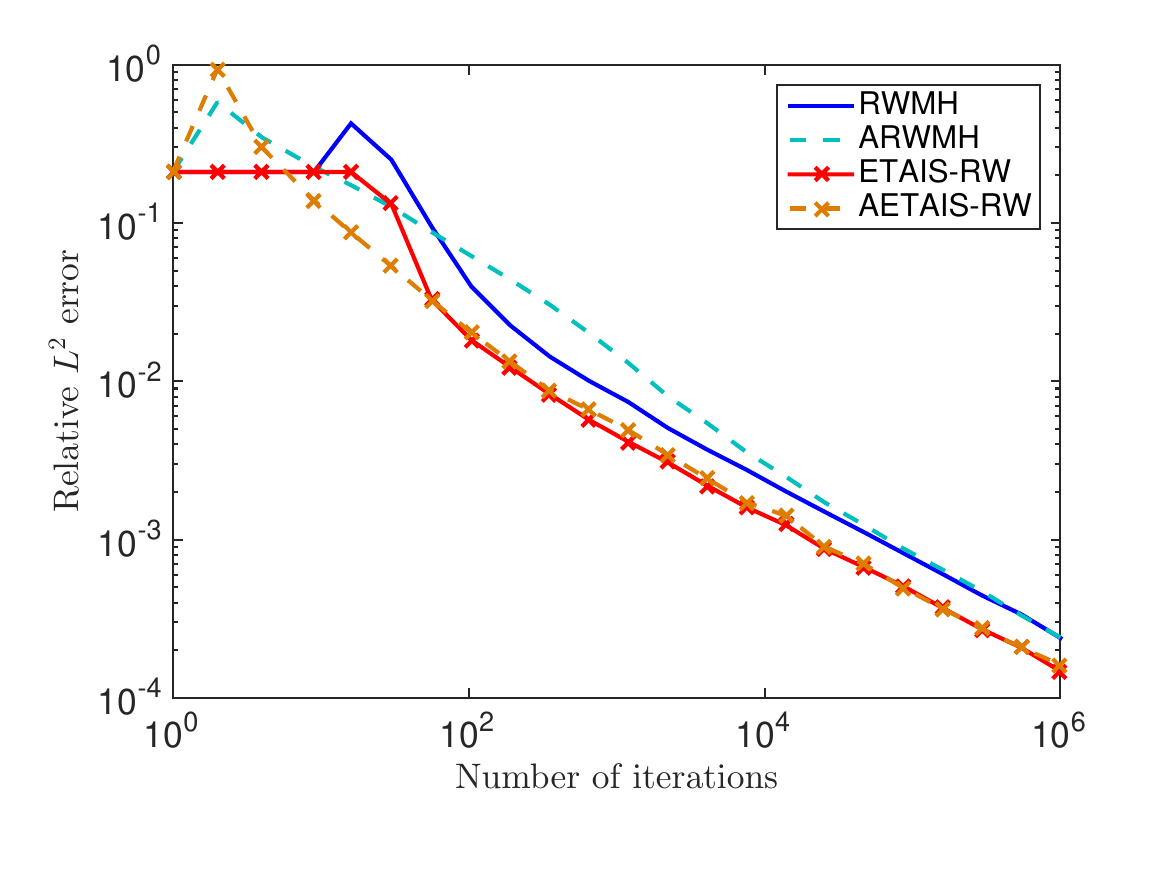}}
\subfigure[Scaling of ETAIS-RW]{\includegraphics[width=0.45\textwidth]{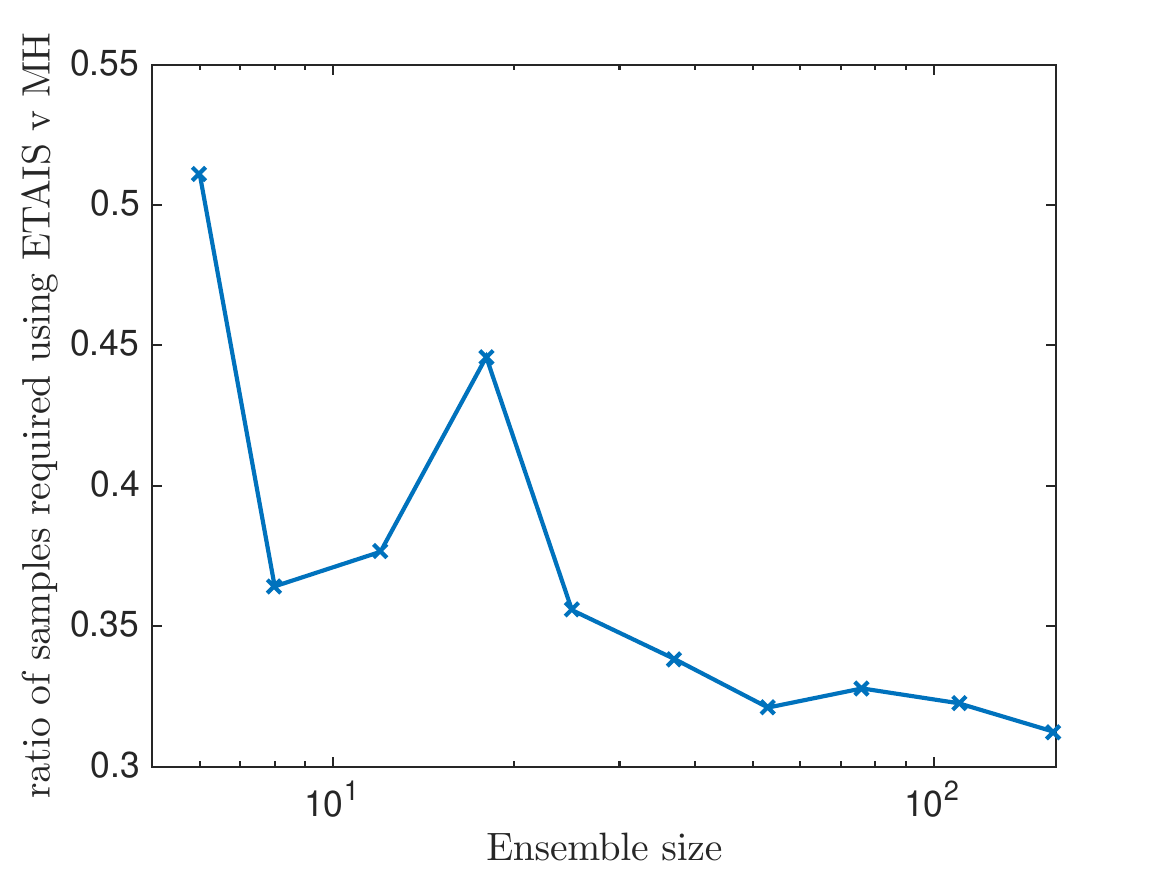}}

\caption{Left: Convergence of the (A)RWMH and (A)ETAIS-RW algorithms
  for a 1D Gaussian. Computed using 50 ensemble members, and the
  ETPF resampler. Right: Ratio of ETAIS-RW samples required to reach the same tolerance
  as the RWMH algorithm for a range of ensemble sizes. }
\label{fig:MH1 L2}
\end{figure}

Figure~\ref{fig:MH1 L2} (left) shows that the ETAIS-RW algorithm converges
to the posterior distribution significantly faster than the RWMH
algorithm, in both $L^2$ error and relative error in the moments. The adaptive algorithms are also shown
in Figure~\ref{fig:MH1 L2}. We see that both adaptive algorithms
converge to the posterior at a similar speed to the respective
optimised algorithm. This shows that, particularly for ETAIS, we can
optimise simulations efficiently on the fly.

Figure~\ref{fig:MH1 L2} (right) was produced by identifying the
optimal value of the scaling parameter for a range of ensemble sizes,
and then computing the ratio of ETAIS samples needed in comparison
with RWMH samples for the same level of error. The decreasing trend shows superlinear improvement of
ETAIS with respect to ensemble size, in terms of the number of
iterations required, which is a demonstration of our belief that
communication between the ensemble members should give us added value over and above that
provided by naive parallelism. This decrease is due to the increasing
effective sample size shown in Figure~\ref{fig:neff-M} (b).

\subsection{Multimodal targets and the effect of resampler quality}
\label{sec:bimodal}
In this second example we investigate the behaviour of the ETAIS
algorithm when applied to a bimodal problem, given by a mixture of two
Gaussians,
\[\pi(x) = 0.2\pi_1(x) + 0.8\pi_2(x),\]
where $\pi_1$ is the density of a $\mathcal{N}\left (\begin{pmatrix} 1
    \\ 1 \end{pmatrix}, \begin{pmatrix} 0.1 & 0 
    \\ 0 & 0.1 \end{pmatrix} \right )$ random variable and $\pi_2$ is the
  density of a $\mathcal{N}\left (\begin{pmatrix} -5
    \\ -5 \end{pmatrix}, \begin{pmatrix} 2.75 & -2.25 
    \\ -2.25 & 2.75 \end{pmatrix} \right )$ random variable.

MH methods can struggle
with multimodal problems, particularly where switches between the
modes are rare, resulting in incorrectly proportioned modes in the
histograms. This example demonstrates that the ETAIS algorithm
redistributes samples to new modes as they are found. This means
that we expect the number of samples in a mode to be approximately
proportional to the probability density in that mode, resulting in
faster global convergence. In particular, we will look at the effect
of using either the ETPF, MT or a standard bootstrap resampler.

\begin{figure}[htb]
\centering 
\subfigure[Target density]{\includegraphics[width=0.45\textwidth]{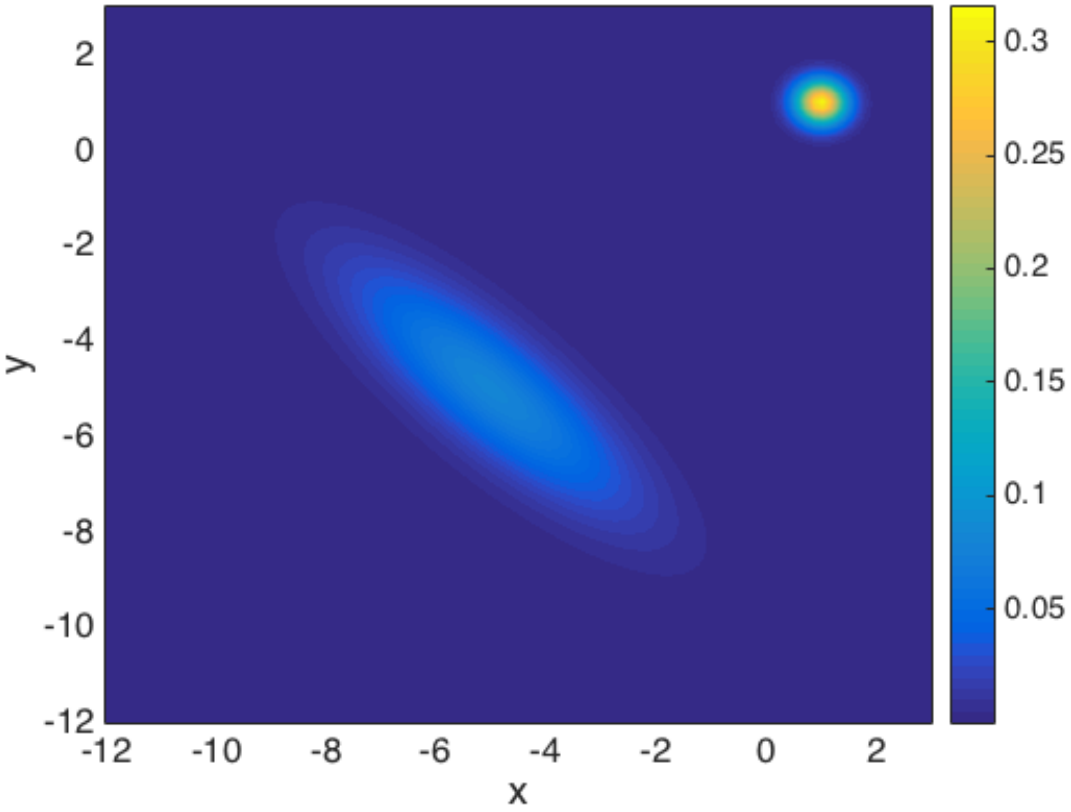}}
\subfigure[Comparison of resamplers]{\includegraphics[width=0.45\textwidth]{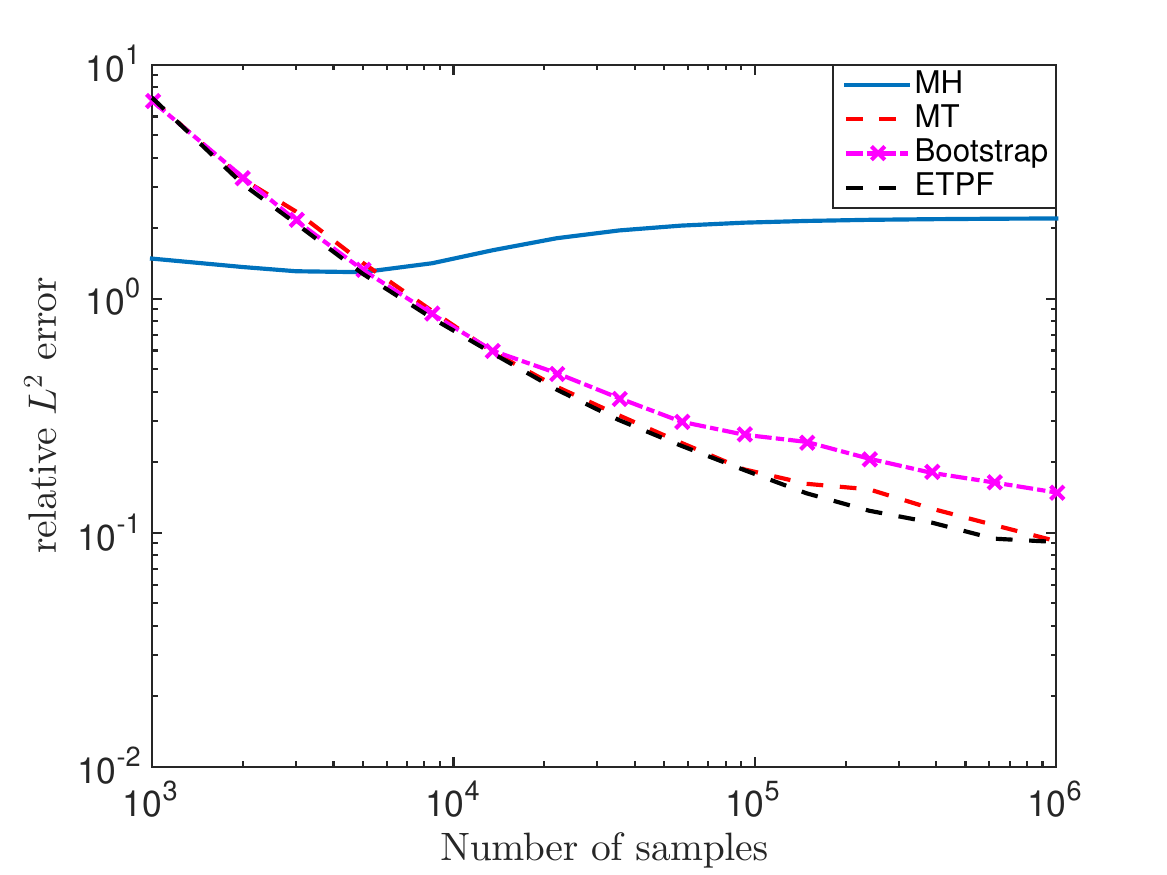}}

\caption{Left: Bimodal target density. Right: Convergence of RWMH and
  ETAIS-RW with ETPF, MT and bootstrap resamplers. }
\label{kshdbf}
\end{figure}

Figure \ref{kshdbf} shows the bimodal target density, and convergence
plots for the RWMH and ETAIS-RW with ETPF, MT and bootstrap
resamplers, averaged over 32 repeats. As expected, a higher quality resampler leads to better
proposal distributions, which in turn leads to greater stability and
more reliable convergence. However, this also shows that the MT is a
good greedy approximation of the ETPF, and for a fraction of the cost
when the number of ensemble members is larger, as shown in Section
\ref{sec:MT}. The RWMH algorithm fails to converge efficiently since chains only very rarely
make switches between the two modes, leading to very slow mixing. This
demonstrates the advantage that ensemble-based methods can have over
other methods when the target
is multimodal.

\subsection{A Mixture Model}\label{sec:mixture}

The technique of mixture modelling employs well known parametric
families to construct an approximating distribution
which may have a complex structure. Most commonly, Gaussian kernels
are used since underlying properties in the data can often be assumed
to follow a Gaussian distribution. An example would be if a
practitioner were to measure the heights of one hundred adults, but
failed to record their gender. The data could be considered as one
population with two sub populations, male and female. The problem then
might be to find the average height of adult males from the data. In
this case, since height is often considered to follow a Gaussian
distribution, it makes sense to model the population as a mixture of
two univariate Gaussian distributions.

A well known problem in the Bayesian treatment of mixture modelling is
that of identification, sometimes referred to as the label-switching
phenomenon. The likelihood distribution for mixture models is
invariant under permutations of the mixture labels. If a mixture has
$n$ means and the point $(\mu_1, \dots, \mu_n)$ maximises the
likelihood, then the likelihood will also be maximised by
$(\mu_{\varphi(1)}, \dots, \mu_{\varphi(n)})$ for all permutations
$\varphi(\cdot)$. This means that the number of modes in the posterior
distribution is of order $\mathcal{O}(n!)$. As we have seen it can be
hard for standard MH algorithms to obtain reliable inference for
posterior distributions with a large number of modes, or even a small
number of modes which are separated by a large distance.

\subsubsection{Target Distribution}\label{sec:mixture_target}

In particular we look at a data set where we assume that there are two
subpopulations within the overall population. Since both
subpopulations will be approximated by Gaussian distributions we have
five parameters which we need to be estimated, two means $\{\mu_{1}, \mu_{2}\}$,
two variances $\{\sigma^2_{1}, \sigma^2_{2}\}$, and the probability, $p$, that an
individual observation belongs to the first subpopulation. We have 100
data points, $D_i$, which we assume to be distributed according to
\[
	D_i \sim p\mathcal{N}(\mu_1, \sigma^2_1) + (1-p)\mathcal{N}
	(\mu_2, \sigma^2_2), \quad i = 1,\dots,100,
\]
where $p \in [0, 1]$, $\mu_{1},\mu_2 \in \mathbb{R}$ and $\sigma^2_{1},\sigma^2_2
\in \mathbb{R}^+$. Due to the domains of these parameters and also
some prior knowledge, we assign the priors
\[
	p \sim \text{Beta}(1,1), \quad \mu_{1,2} \sim \mathcal{N}(0, 4)
	 \quad \text{and} \quad \sigma^2_{1,2} \sim \text{Gamma}
	(\alpha=2, \beta=1).
\]
If we collect these parameters in the vector $\theta = (p, \mu_1,
\sigma^2_1, \mu_2, \sigma^2_2)^\top$, the resulting posterior
distribution is
\begin{equation}\label{eq:mixture_posterior}
	\pi(\theta|D) \propto \prod\limits_{i=1}^{100} (p\mathcal{N}
	(D_i;\mu_1, \sigma^2_1) + (1-p)\mathcal{N}(D_i;\mu_2, \sigma^2_2))
	\prod\limits_{i=1}^5 \pi_0^i(\theta_i),
\end{equation}
where $\pi_0^i(\cdot)$ is the prior density function corresponding to $\theta_i$.

\begin{figure}[htb]
\centerfloat
\includegraphics[width=1.2\textwidth]{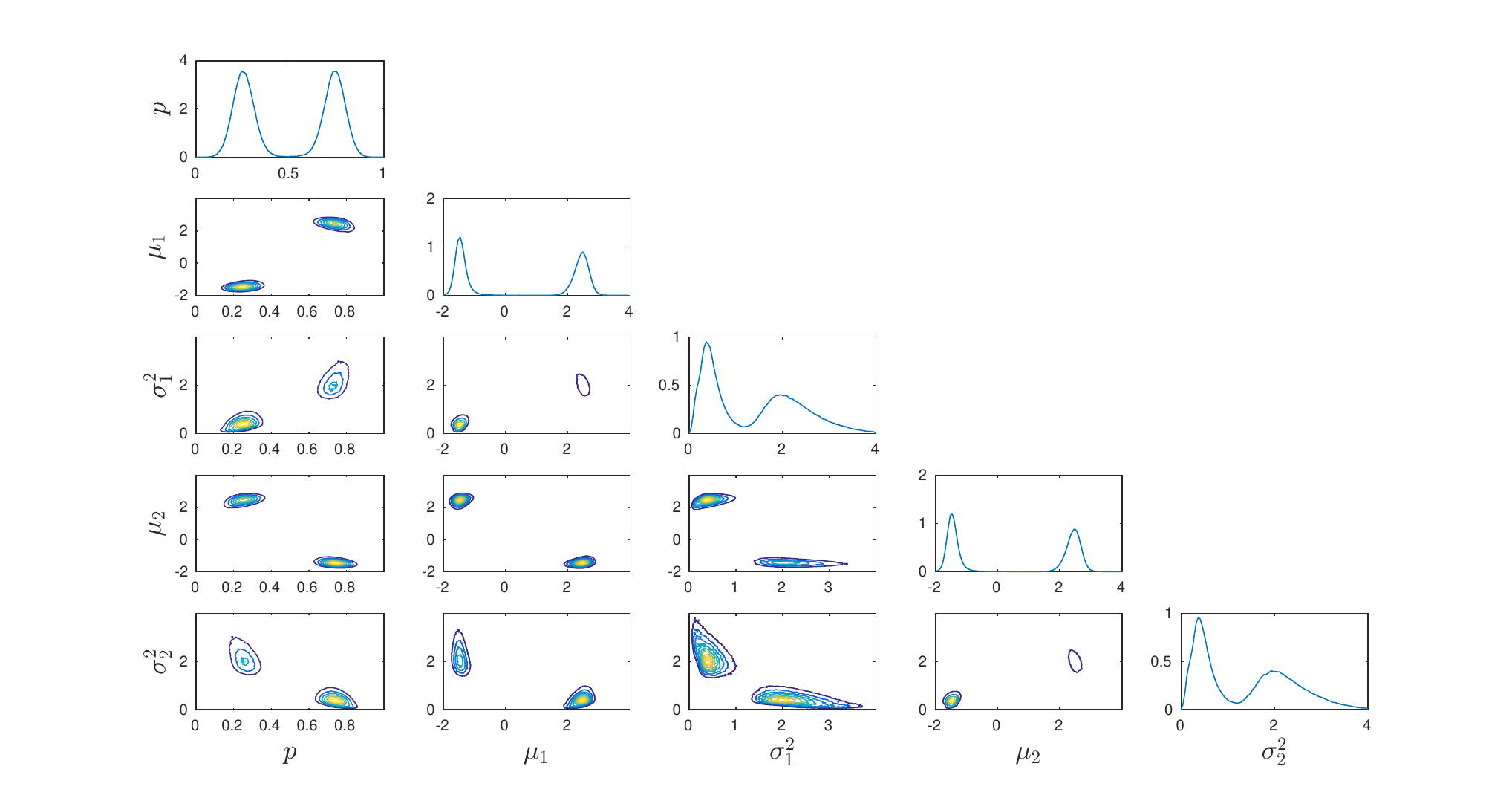}
\caption{The posterior distribution of $\theta$ as given in
Equation~\eqref{eq:mixture_posterior} as found from 10 million samples
from the ETAIS algorithm. The main diagonal contains the marginal
distributions of each $\theta_i$, and the lower triangular contours
represents the correlation between pairs of parameters.}
\label{fig:mixture_posterior}
\end{figure}

Figure \ref{fig:mixture_posterior} presents a visualisation of the
posterior distribution for this problem, created from 10 million
samples produced by the ETAIS-RW algorithm.

%%%%%%

\subsubsection{Implementation}
\label{sec:mixture_implementation}
In order to quantify error, we note that the probability density should
be evenly divided between the two modes. This is due to the symmetric
prior for $p$, and and the same priors being assigned to $\mu_1$ and
$\mu_2$, and also to $\sigma^2_1$ and $\sigma^2_2$. To decide which
mode a sample belongs to we define a plane which bisects the posterior
so that each point on this plane lies exactly halfway between the two
true solutions to the inverse problem i.e. the value of $\theta$ used
to generate the data, and also the $\theta$ obtained by a relabelling
of the parameters. Now that we can assign a sample to a particular
mode, we can calculate the density in each mode by summing the weights
associated to all samples in that mode,
\[
	\bar{w}_k = \sum\limits_{i=1}^N w_iI_{X_i \in \text{Mode $k$}}, \quad k = 1, 2,
\]
and the relative error in the amount of density in each mode is then
\begin{equation}\label{eqn:mode_prop}
	w_\text{error} = 2\left|\frac{\bar{w}_1}{\bar{w}_1+\bar{w}_2} - \frac{1}{2}\right|.
\end{equation}

Since the probability $p$ is constrained to lie in the interval
$[0,1]$, and the variances must be positive, it can be wasteful to use
Gaussian proposal distributions, which will produce samples outside of
the support of the posterior. Moreover, the value of the variances of
the unknown distributions are strictly non-negative. The algorithm
will be most efficient if the proposal and posterior distributions are
mutually absolutely continuous. It is also useful to be able to pick
proposal densities for which the variance is easily scaled so that we
can tune them to optimise efficiency. Thus we pick the following
proposal distributions for the $p$, $\mu$ and $\sigma^2$ parameters, respectively;
\[
	q_p \sim \text{Beta}(\delta^{-2}p, \delta^{-2}(1-p)), \quad q_{\mu_{1,2}} \sim \mathcal{N}(\mu_{1,2}, 4\delta^2) \quad \text{and} \quad q_{\sigma^2_{1,2}} \sim \text{Gamma}(\alpha^*, \beta^*),
\]
where $\alpha^* = \sigma^2_{1,2}\beta^*$, $\beta^* =
\sigma^2_{1,2}/2\delta^2$ and $\delta$ is a scaling parameter to be
tuned. This means that our proposal distributions will not be a
mixture of multivariate Gaussians, but independent mixtures of
univariate Beta, Gamma and Gaussian distributions.

In the numerics which follow we have increased the ensemble size from
$M=50$ to $M=500$ to compensate for the increase in dimension. We also
use the MT algorithm for the resampling step because of the reduced
computational cost. The method otherwise remains the same as in
previous examples. We perform test runs to find the optimal scaling
parameters considering convergence to modes with equal density. We
then calculate the convergence rates of the algorithms by producing 10
million samples from the posterior with each algorithm, and repeat the
simulation 32 times.

% numerics.
\subsubsection{Convergence of MH vs ETAIS}\label{sec:mixture_conv}

\begin{table}[!htb]
      \centering
        \begin{tabular}{|l|r|r|}
	\hline
	Algorithm	& MH & ETAIS \\ \hline
	$\delta^*$	& 1.3e-1     & 2.3e-1 \\
	\hline
	\end{tabular}
	\vspace{1mm}
	\caption{Optimal values of the scaling parameter. The MH algorithm is optimised using the acceptance rate, and the ETAIS algorithm is optimised using the effective sample size.}
	\label{table:mixture_opt_beta}
\end{table}

The optimal scaling parameters for the MH and ETAIS algorithms with the proposal distributions described in Section~\ref{sec:mixture_implementation} are given in Table~\ref{table:mixture_opt_beta}.

\begin{figure}[htb]
\centering
\includegraphics[width=0.8\textwidth]{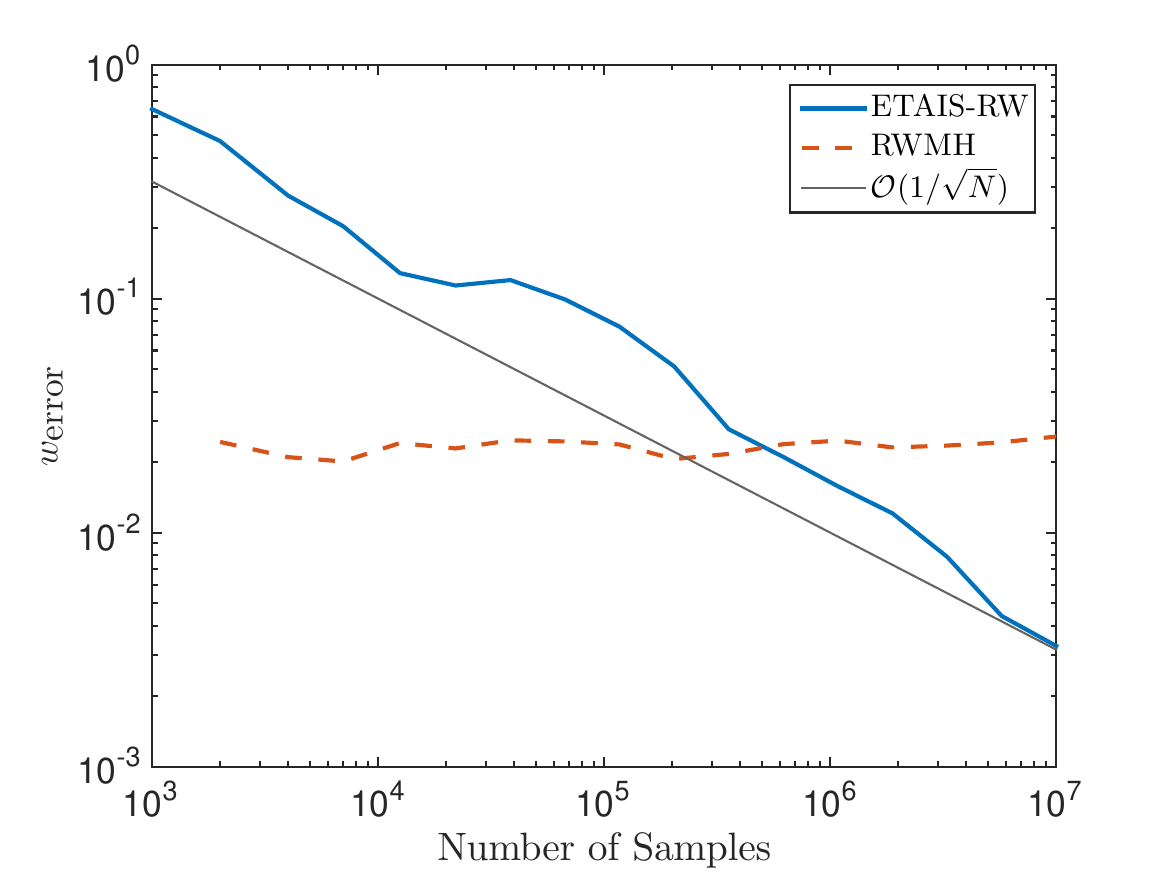}
\caption{Convergence of the ETAIS algorithm for the mixture model described in Section~\ref{sec:mixture}, convergence judged using the criterion in Equation~\eqref{eqn:mode_prop}. Implementation described in Sections~\ref{sec:mixture_implementation} and \ref{sec:mixture_conv}. Resampling is performed using the MT scheme.}
\label{fig:mixture_modes}
\end{figure}

Convergence of the relative error for the two algorithms is displayed
in Figure~\ref{fig:mixture_modes}. ETAIS converges at the expected
$\mathcal{O}(1/\sqrt{N})$ rate, whereas the MH algorithm converges to
locally smooth histograms but with the wrong proportion of samples in
each mode. The relatively low value of the error for the MH example is
due to the priors covering the sample space evenly, however since
transitions are near impossible with a small value of the scaling
parameter, this error will take a very long time to reduce. This
problem was discussed in Section~\ref{sec:bimodal}.

In this example, we have only considered a relatively low dimensional
mixture model problem, with only two distributions in the
mixture. With more elements in the mixture, and/or an undefined number
in the mixture, the dimension of this problem will quickly
increase. Importance sampling schemes such as this suffer from the
curse of dimensionality, limiting the size of the dimension of target
density that it can efficiently sample from. However, this example
does demonstrate the remarkable convergence of the ETAIS for multimodal
target distributions.

\subsection{Data assimilation with Lorenz `63 trajectories}\label{sec:lorenz}

The sampling algorithm we have introduced in this paper is designed for problems where the likelihood density function is expensive to calculate, and the computational overhead involved in the calculation of the importance weights is dwarfed. This situation commonly occurs in inverse problems where the model being investigated involves the solution of a differential equation. In this example we attempt to recover the initial position of a particle with motion governed by Lorenz's 1963 atmospheric convection equations~\cite{lorenz1963deterministic}:
\begin{align*}
	\frac{\text{d}x}{\text{d}t} &= \sigma(y-x), \\
	\frac{\text{d}y}{\text{d}t} &= x(\rho-z)-y, \\
	\frac{\text{d}z}{\text{d}t} &= xy - \beta z.
\end{align*}
When the parameters are chosen to be $\rho = 28$, $\sigma=10$ and $\beta=8/3$ this system has chaotic solutions. If we are interesting in calculating the initial position by observing the location of the particle at certain points in time, it can quickly become intractable when there is any noise in the observations.

\begin{figure}[htb]
\centering
\includegraphics[width=0.8\textwidth]{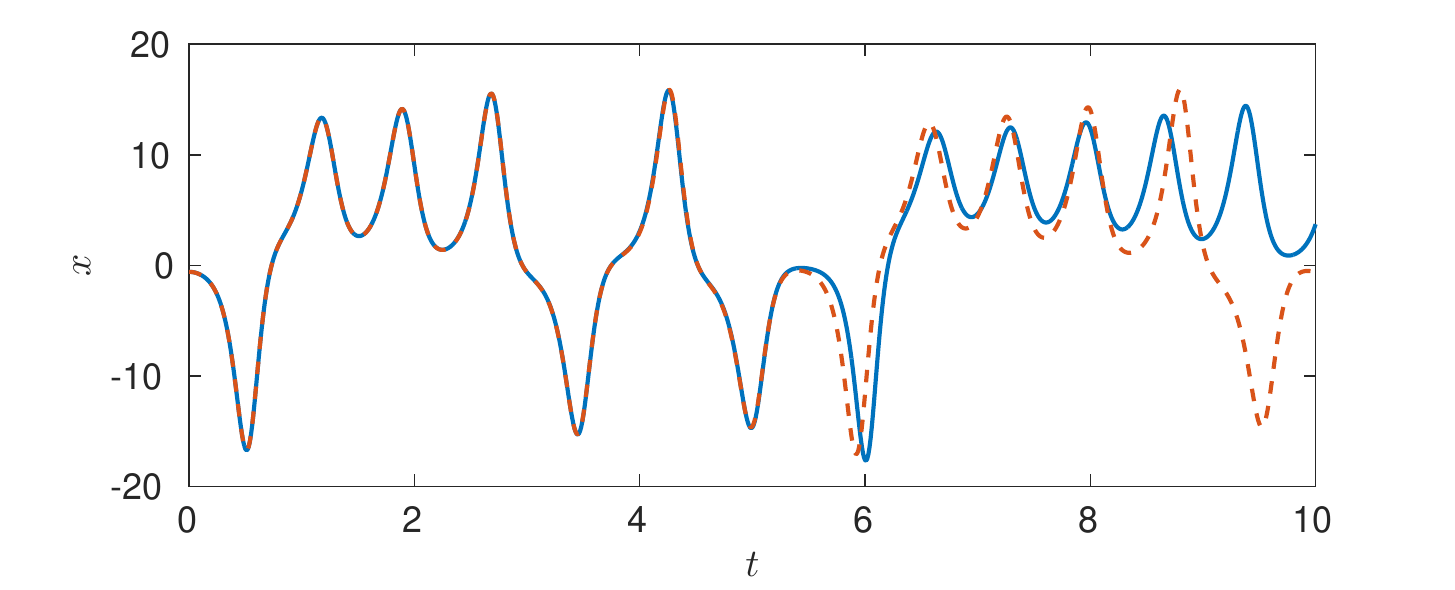}
\caption{Trajectory of the $x$ component of a particles position varying in time when motion is governed by the Lorenz `63 equations. The two trajectories are calculated using an initial position of $[-0.587,-0.563,16.870]$ and $[-0.590,-0.563,16.870]$.}
\label{fig:Lorenz_xtraj}
\end{figure}

To demonstrate how small errors in the initial condition of a Lorenz solution can affect the trajectory of a particle, Figure~\ref{fig:Lorenz_xtraj} shows the path taken by two particles with very similar initial conditions,
\[
	\mathbf{x}_0^1 = \begin{pmatrix} -0.587\\-0.563\\16.870 \end{pmatrix} \quad \text{and} \quad \mathbf{x}_0^2 = \begin{pmatrix} -0.590\\-0.563\\16.870 \end{pmatrix}.
\]
This small difference in the $x$-dimension of the initial condition leads to the trajectories of the two particles decoupling near $t=5$. This chaotic behaviour means that if we take only a few observations over a long time period we will find that many trajectories which may be vastly different achieve similar values of the likelihood function and so the posterior distribution for the initial condition can become very complex.%

\subsubsection{Target Distribution}\label{sec:lorenz_target}

For this example, we observe noisily the position of a particle at ten
equally spaced points in the time interval $t \in (0, 1]$. The Lorenz
equations are solved with the chaotic parameters given above and the
initial condition $\mathbf{x}_0^1$. A time step of $h=1\times 10^{-3}$
was used to evolve the equations numerically with the explicit Euler
method. The noise added to each observation was taken from the
distribution $\mathcal{N}(0, 0.1^2I)$. Priors for the initial condition
coordinates were given by
\[
	x_0 \sim \mathcal{N}(-0.5, 0.4^2), \quad y_0 \sim\mathcal{N}(-0.5, 0.4^2) \quad \text{and} \quad z_0 \sim \mathcal{N}(15, 0.4^2).
\]
The posterior density function takes the form
\[
	\pi(\mathbf{x}_0|D) \propto \exp\left\{-\frac{1}{2}\|\mathcal{G}(\mathbf{x}_0) - D\|^2_\Sigma\right\}\mathcal{N}(\mathbf{x}_0; m, 0.4^2\text{I}),
\]
where $\mathcal{G}$ calculates the solution to the Lorenz equations, using the Euler method with a time step of $h=1\times 10^{-3}$ starting at the initial condition $\mathbf{x}_0$, and returns the position at $t = 0.1, ..., 1$. The mean of the prior distribution $m = [-0.5, -0.5, 15]$.

\begin{figure}[htb]
\centering
\includegraphics[width=0.8\textwidth]{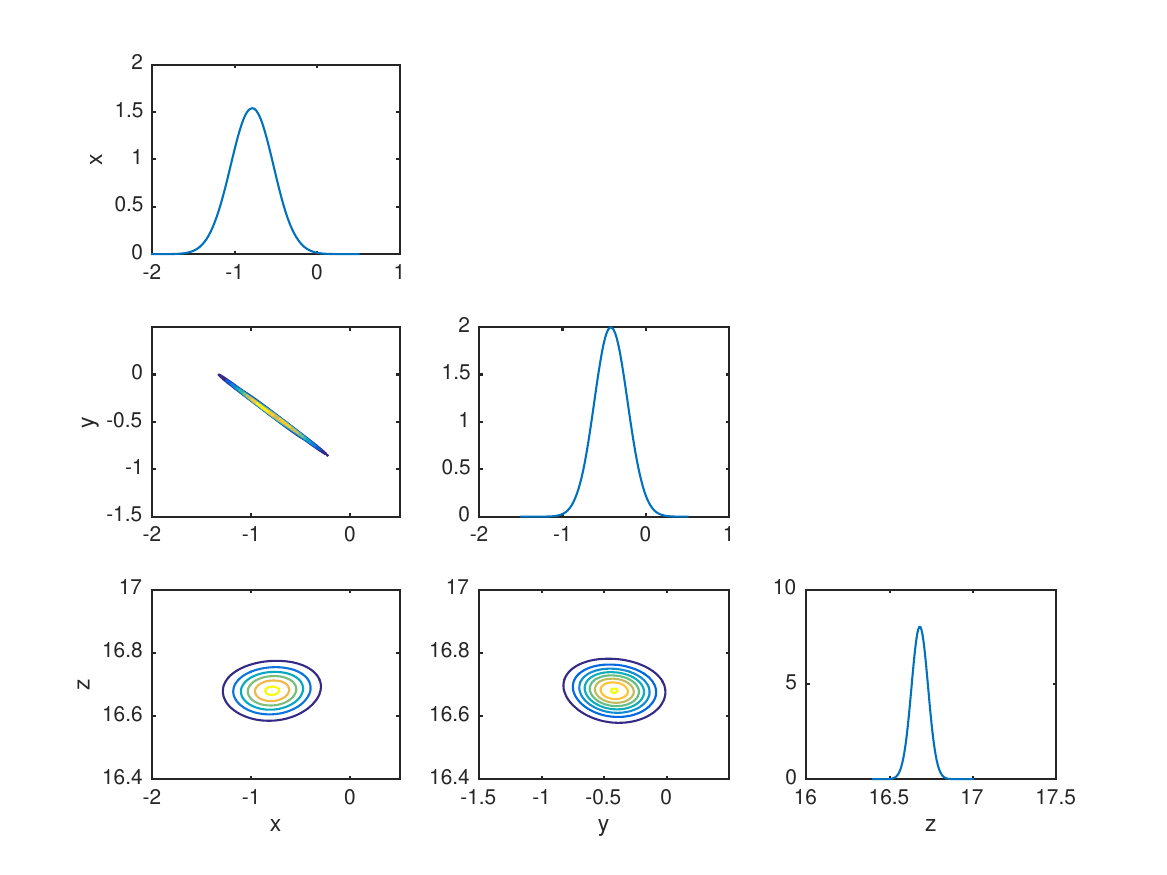}
\caption{Posterior distribution for the initial condition of the particle with motion governed by the Lorenz equations as described in Section~\ref{sec:lorenz_target}.}
\label{fig:Lorenz_posterior}
\end{figure}

Figure~\ref{fig:Lorenz_posterior} shows the marginal distributions of
the posterior distribution for the initial condition $\mathbf{x}_0$ as
found using the MH algorithm. While the $z$-dimension is largely
uncorrelated with the other two dimensions, the correlation between
$x$ and $y$ is -0.97, which makes it a challenging posterior distribution to explore without a transformation of the parameter space.

%%%%%%
\subsubsection{Implementation}\label{sec:lorenz_implementation}

As in the previous example, we have no analytic form for the normalisation constant of the posterior and so we will measure convergence of the algorithms by calculating convergence to the posterior mean, where the truth is calculated using a very long chain produced by the MH algorithm.

For this example, we use an ensemble size of $M=1500$, and the MT resampling algorithm. The optimal scaling parameters are calculated for both algorithms using test runs which are not included in the convergence cost calculations. Convergence graphs are produced using 32 repeats of each algorithm with each simulation producing one million samples.

% numerics.
\subsubsection{Convergence of MH vs ETAIS}\label{sec:lorenz_conv}

% The optimal scaling parameters for the MH and ETAIS algorithms are given in Table~\ref{table:Lorenz_opt_parameters}.

% \begin{table}[!htb]
%       \centering
%         \begin{tabular}{|l|r|r|}
% 	\hline
% 	Algorithm	& MH & ETAIS \\ \hline
% 	$\delta^*$	& 1.0e-1     & 1.3e-2 \\
% 	\hline
% 	\end{tabular}
% 	\vspace{1mm}
% 	\caption{Optimal values of the scaling parameter. The MH algorithm is optimised using the acceptance rate, and the ETAIS algorithm is optimised using the effective sample size.}
% 	\label{table:Lorenz_opt_parameters}
% \end{table}

\begin{figure}[htb]
\centering
\includegraphics[width=0.7\textwidth]{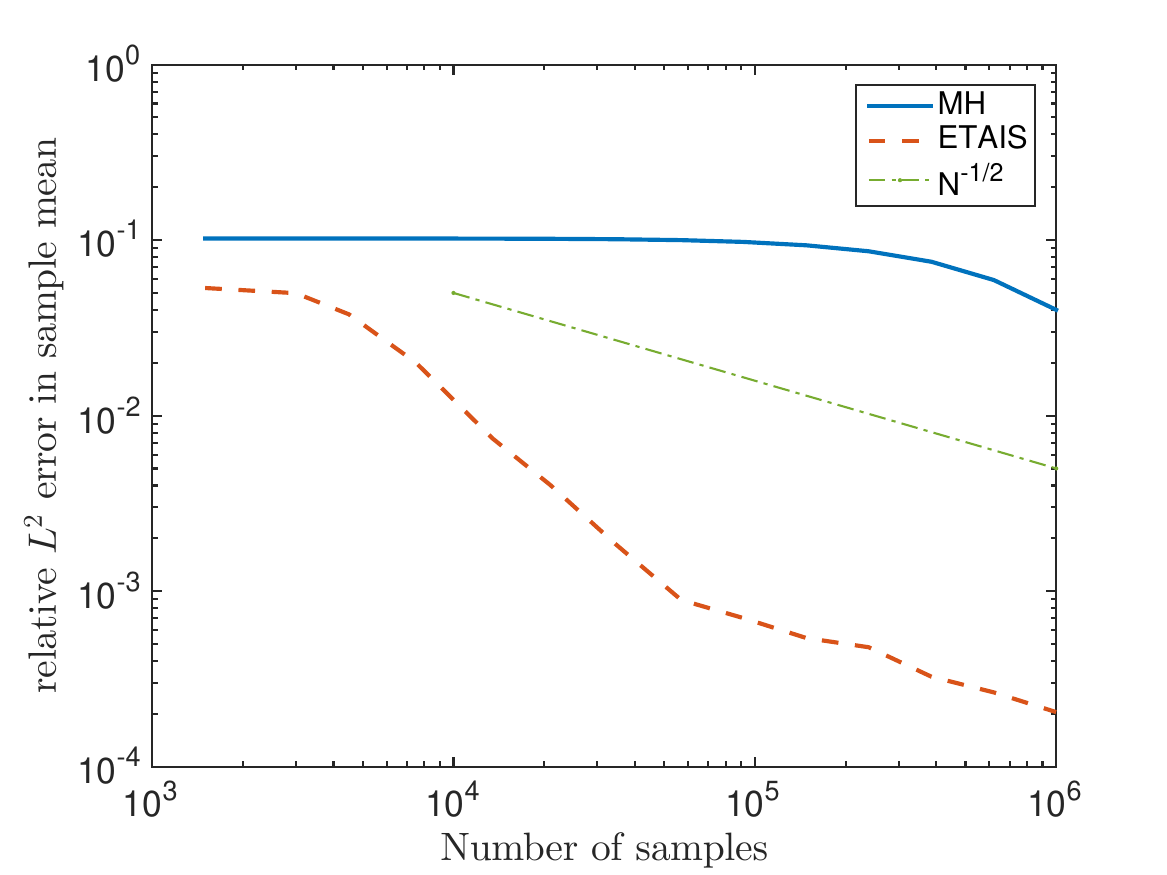}
\caption{Convergence of the sample mean to the posterior mean for the MH and ETAIS algorithms. Posterior distribution as described in Section~\ref{sec:lorenz_target}.}
\label{fig:Lorenz_convergence}
\end{figure}

Figure \ref{fig:Lorenz_convergence} shows convergence plots for
  this problem using RWMH and ETAIS. The burn-in of the ETAIS is
  remarkably better than MH, as the ensemble quickly covers the
  manifold on which the majority of the density lies close to. The
  plot shows convergence with respect to number of samples, but the
  two approaches have differing costs, due to the calculation of the
  denominator in the weights in the ETAIS, and the resampling step
  (MT). Assuming that both error curves decay like $N^{-1/2}$ from
  $N=10^6$, then for a given error tolerance, the number $N_{\rm MH}$
  of MH samples required is given by
\[N_{\rm MH} = N_{\rm ETAIS} \times 3.838 \times 10^4,\]
where $N_{\rm ETAIS}$ is the number of ETAIS samples required for the
same error tolerance. The cost per sample for the MH and ETAIS methods
for this problem are $1.64 \times 10^{-5}$ and $1.46 \times 10^{-3}$
respectively. Here the ETAIS results were produced with a serial
implementation, so this cost may be slightly overoptimistic, since the
runtimes do not include extra overheads of communicating the ensemble
states between processors. However, putting this together, we arrive at a speed-up factor of
over 430, which will undoubtedly overshadow any such underestimate of
the cost-per-sample of the ETAIS approach. This demonstrates the benefits of this approach for
problems with expensive likelihoods and challenging posteriors such as
this one.

\section{Discussion and Conclusions}\label{Sec:Conc} 

We have explored the application of 
low dimensional Bayesian inverse problems. We have demonstrated numerically
that this method converges faster than the analogous naively parallelised
Metropolis-Hastings algorithms. Further experimentation with the Metropolis
Adjusted Langevin Algorithm (MALA), preconditioned Crank-Nicolson (pCN),
preconditioned Crank-Nicolson Langevin (pCNL) and Hamiltonian
Monte Carlo (HMC) proposals has yielded similar results\cite{Paul}.

The implementations of the ETAIS that have been presented in this
  paper, have been run in serial, but we argue that, particularly for
  inverse problems with expensive likelihoods, for example in the case
  where this requires the approximation of the solution of a
  differential equation, that this method is an excellent candidate
  for parallelisation, with the likelihood evaluations for each
  particle being computed on a different processor. This said, the numerical results that we have presented demonstrate
  that even a serial implementation can outperform standard
  Metropolis-Hastings methods.
  
% Importantly, we have compared the efficiency of our scheme
% with a naive parallelisation of serial methods. Thus our increase in
% efficiency is over and above an $M$-fold increase, where $M$ is the
% number of ensemble members used. Our approach
% demonstrates a better-than-linear speed-up with the number of ensemble
% members used. 

The ETAIS has a number of favourable features, for example the
algorithm's ability to redistribute, through the resampling regime,
the ensemble members to regions which require more exploration. This allows the
method to be used to sample from complex multimodal distributions.

Another strength of the ETAIS is that it can also be used with any MCMC
proposal. There are a growing number of increasing sophisticated MCMC
algorithms (non-reversible Langevin/HMC proposals, Riemann manifold MCMC etc) which could be
incorporated into this framework, leading to even more efficient
algorithms, and this is another opportunity for future work. 

One limitation of the ETAIS approach as described above is that a
direct solver of the ETPF problem (such as FastEMD \cite{FastEMD}) has computational cost
$\mathcal{O}(M^3\log M)$, where $M$ is the number of particles in the
ensemble. As such, we introduced a more approximate resampler the
approximate multinomial resampler, which allows us to push the
approach to the limit
with much larger ensemble sizes. The ETAIS framework is very flexible
in terms of being able to use any combination of proposal
distributions and resampling algorithms that one wishes.

We have demonstrated that the framework that we have considered, with
the use of state-of-the-art optimal transport-based resampling, can
reduce the number of likelihood evaluations required to characterise
complex posterior distributions in low dimensions to a given
degree. We have also introduced a greedy approximation to this
resampler, which drastically reduces the cost, at the loss of some
accuracy, which can be mitigated with the use of larger
ensembles. We have detailed how scaling parameters in the MCMC
proposals that are used within the mixture distribution can be quickly
and efficiently
tuned in an automated way. Lastly we have demonstrated that when the
likelihood is expensive, for instance because it involves the
numerical approximation of the solution of a differential equation, 
that we can achieve orders of magnitude reductions in cost to reach a
given error tolerance in comparison with standard Metropolis-Hastings approaches.

However, tuning the variances and
covariances of the mixture proposal components
globally is likely to be of limited use for multimodal problems where
the modes have very different covariances, or where there are curved
ridges in the density. With increased dimensionality and/or
  complexity of the target the distribution, we also require increases
  in the size of the ensemble if we wish to have a stable ETAIS
  implementation. After a certain point, the extra overheads
  associated with a very large ensemble will outweigh the advantages
  of this approach. However, there are a variety of possible solutions to
this problem that would be worthy of future consideration, not least
the potential to use transport maps\cite{el2012bayesian,parno2014transport} to map Gaussian mixtures to
highly complex non-Gaussian approximations of the posterior. Such a
map could encode local covariance information, and lead to accelerated
and stable ETAIS sampling with much smaller ensemble sizes.

\begin{appendix}

\section{Glossary of acronyms}

\begin{table}[h!]
\centering
\begin{tabular}{|l|l|}
\hline
Acronym & Full name                              \\ \hline
MCMC    & Markov chain Monte Carlo               \\
RWMH    & Random walk Metropolis-Hastings        \\
MALA    & Metropolis-adjusted Langevin algorithm \\
AIS     & Adaptive importance sampling           \\
PMC     & Population Monte Carlo                 \\
ETAIS    & Ensemble transport adaptive importance sampling  \\
ETPF    & Ensemble transport particle filter     \\
MT     & Multinomial transformation      \\
ETAIS-X  & ETAIS with X kernels                    \\
AETAIS-X & Adaptive ETAIS with X kernels          \\ \hline
\end{tabular}
\end{table}

\end{appendix}

\bibliographystyle{siam}
\bibliography{refs}

\end{document}